\documentclass[12pt,a4paper,dvipsnames]{article}

\usepackage[margin=1in]{geometry}

\usepackage[utf8]{inputenc}

\usepackage{amsmath}
\usepackage{amsfonts}
\usepackage{amssymb}
\usepackage{graphicx}
\usepackage{tabularx}
\usepackage{xspace}
\usepackage{booktabs}
\usepackage{hyperref}
\usepackage{multirow}
\usepackage{subcaption}

\usepackage{rotating}

\usepackage[square,numbers]{natbib} 

\usepackage{setspace}

\allowdisplaybreaks

\newcommand{\Nesp}{SetN\xspace}
\newcommand{\NespOne}{SetNLow\xspace}
\newcommand{\NespTwo}{SetNHigh\xspace}
\newcommand{\Brill}{SetB\xspace}
\newcommand{\BrillOne}{SetBLow\xspace}
\newcommand{\BrillTwo}{SetBHigh\xspace}
\newcommand{\Tric}{SetT\xspace}

\begin{document}
\title{Planning profitable tours for field sales forces: A unified view on sales analytics and mathematical optimization}
\author{Anne Meyer\footnote{TU Dortmund University, Leonhard-Euler-Straße 5, D-44227 Dortmund, anne2.meyer@tu-dortmund.de}, Katharina Glock\footnote{FZI Research Center for Information Technology, Haid-und-Neu-Str. 10-14, D-76131 Karlsruhe, kglock@fzi.de}, 
Frank Radaschewski\footnote{PTV Group, Haid-und-Neu-Str. 15, D-76131 Karlsruhe, frank.radaschewski@ptvgroup.com}}
\date{}

\maketitle

\begin{abstract}
{A key task of sales representatives in operational planning is to select the most promising customers to visit within the next days.}
{A strongly varying set of scoring methods predicting or approximating the expected response exists for this customer selection phase.} 
However, in the case of field sales forces, the final customer selection is strongly interrelated to the tour planning decisions.
To this end, we formalize variants of the profitable sales representatives tour problem as a multi-period team orienteering problem, thereby providing a unified view on the customer scoring and the tour planning phase.
In an extensive computational study on real-world instances from the retail industry, we systematically examine the impact of the aggregation level and the content of information provided by a scoring method and the sensitivity of the proposed models concerning prediction errors.
We show that the selection of a customer scoring and tour planning variant depends on the {available} data. {Furthermore, we work out where to put effort in the data acquisition and scoring phase to get better operational tours.}
\\

\end{abstract}

\onehalfspacing

\section{Introduction}

In recent years, marketing activities have become ``more individual customer-centered and relationship-driven instead of mass market-centered and transaction-driven'' \citep{Rust2014}.
This holds not only for companies selling directly to end-users (B2C) but also for ``relationship intensive'' \citep{Rust2014} firms selling to other companies (B2B).

As a consequence, the role of sales representatives, visiting customers on-site, changed considerably without losing its important role:
In a study of the consulting company McKinsey, 74\% out of 1'000 large organizations buying goods or services from B2B sellers want to talk to a sales representative if they are interested in a new product or service\footnote{https://www.mckinsey.com/business-functions/marketing-and-sales/our-insights/what-the-future-science-of-b2b-sales-growth-looks-like (12.08.2019)}.
Regular visits for simply taking orders, historically an important task for sales representatives, lose their significance, {while} visits targeted to inform about current promotions, to explain new products or complex services, and to gather knowledge about the needs of customers by maintaining a long-term relationship gain in importance.

In cooperation with PTV Group, an IT provider of optimization services for sales force management, {in this work we tackled} the operational planning of sales forces visiting customers on-site in the retail and consumer goods industry.
Both in literature and in practice, this task has been approached from two distinct perspectives {and in completely different research communities}:
On the one hand, there are optimization algorithms for determining daily visit sequences for sales representatives considering restrictions such as working times {or visit time windows}.
On the other hand, there is a wide range of sales analytics techniques for assessing the importance of individual customers or customer groups,
{reaching from the simple ABC segmentation of customers based on their historic sales to sophisticated methods for predicting future sales based on extensive databases. 
These methods require vastly different efforts for data acquisition, maintenance, and analysis.
}



{These possibly considerable expenses for sales departments raise the question of which effort is justified by improved operational decisions.}
{This work argues that the answer to the question can only be found by jointly considering the customer scoring (or prediction) phase and the tour planning (or decision) phase.}
Most important {in this context} is the interplay of sales analytics and optimization approaches: 
Without a sufficient assessment of customers' purchasing potential, the operative customer selection is likely to be ineffective.
Inversely, day-to-day operations are subject to various restrictions that limit the impact of the predicted customer contributions on the planned tours and, therefore, the benefit of using expensive data acquisition and prediction methods.

{In this work, we study the combination of sales analytics and tour planning in an extensive case study. To our knowledge, this is the first to address both phases based on real-world data sets from two large B2B sales operations. 
In detail, our work makes the following contributions:}
{(1) We propose a new categorization of how customer scoring methods approximate the customer's profitability.}
(2) We formalize the problem of planning profitable operative tours for sales representatives, i.e. tours for the next couple of days, with a special focus on a unified view on {categories of} customer scoring methods and {corresponding} tour planning models.
(3) Based on a comprehensive computational study on instances from the retail industry, we investigate the interplay between expressing profitability of a customer and the tour planning task {by synthetically generating data for all profitability approximation categories.} 
We also investigate the sensitivity of tour planning models considering errors in the customer scoring stage.
{For solving the tour planning problems we apply a state-of-the-art framework for team orienteering problems.}
(4) Furthermore, we discuss important insights for selecting customer scoring and tour planning models that fit the data situation of a company, and we derive findings on where to focus the effort in the prediction phase to enable better operational tours.

The remainder of this paper is organized as follows:
In Section \ref{sec:problem_def}, we introduce the profitable tour planning problem for sales representatives. 
Section \ref{sec:related_work} summarizes the related work in the fields of optimization and sales analytics. 
We then formally introduce models combining these two aspects in Section \ref{sec:model}.
Section \ref{sec:instances} summarizes the test instances from industry that are at the basis of the computational study.
Section \ref{sec:ev2mls} gives insights into the performance of the solution algorithms, while the proposed models are studied in detail in Section \ref{sec:evaluation_models}. 
We conclude this paper with key insights in Section \ref{sec:outlook}.

\section{Profitable tours for sales representatives}\label{sec:problem_def}

\begin{figure}[h!]
 \includegraphics[width=0.7\textwidth]{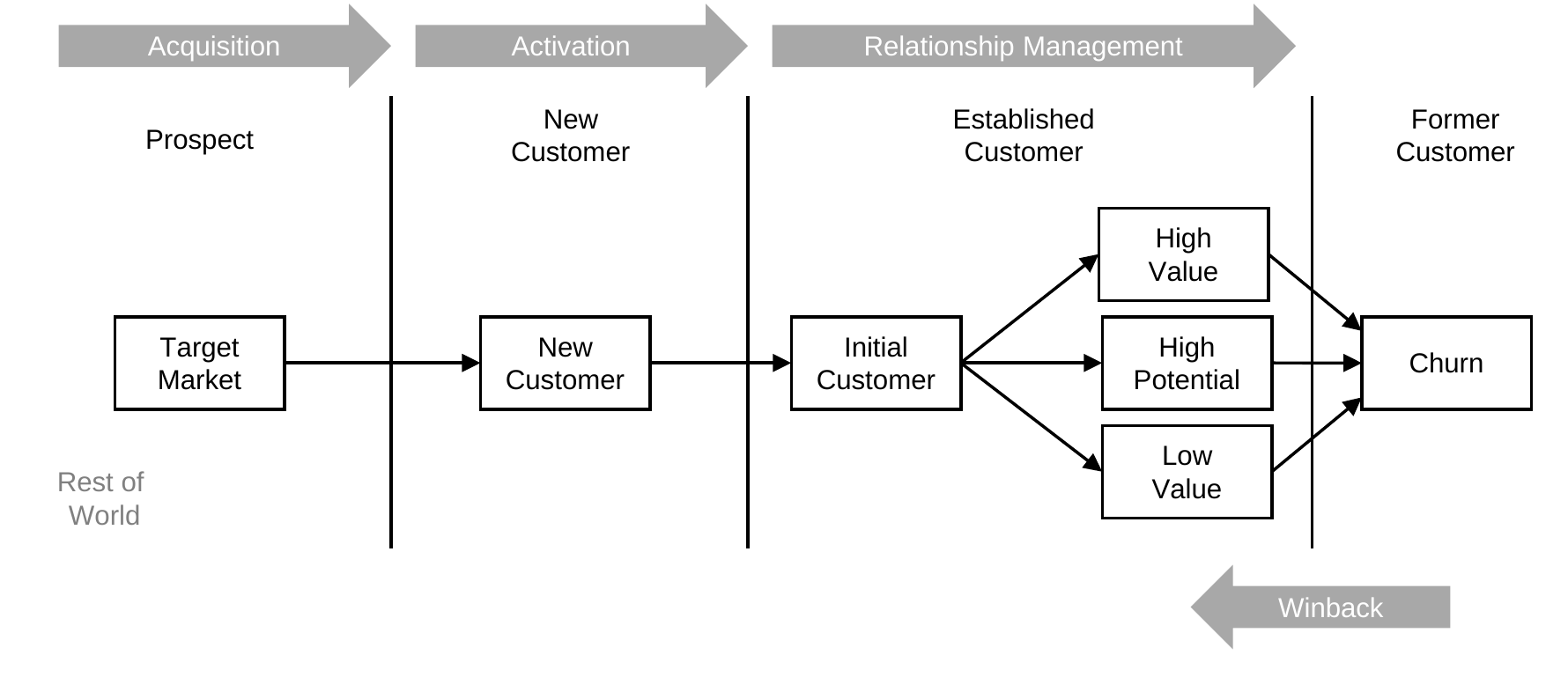} 
  \centering
  \caption{Customer life-cycle containing different stages based on \cite[][p. 29]{Linoff2011}.}
  \label{fig:CustLife}
\end{figure}

The responsibilities of a sales representative can be well explained in the view of a customer life-cycle model such as the one depicted in Figure~\ref{fig:CustLife}:
In general, the customer life-cycle differentiates between prospects and customers.
Following the definitions of \citet[][p. 29]{Linoff2011}, prospects are companies in the target market which are not yet customers, but have shown a first interest, for example by responding to an e-mail.
New customers are companies which have made the first purchase or who have signed a contract, while established customers are customers who return, and for whom the relationship hopefully becomes broader or deeper.
Former customers have left, for example to a competitor\footnote{Note that there exist different customer life-cycle models in literature with a deviating use of terms, e.g., \citet[][p. 44]{Hall2017} or \citet[][p. 53]{Kumar2018}}.

%

Among the clients of PTV group there are companies that operate sales forces with several hundred of sales representatives serving tens of thousands of customers (e.g., supermarkets) in a country such as Germany.
Due to the high costs associated with such a sales force, the key goal of each company is to effectively and efficiently manage its sales representatives.

\subsection{Planning profitable tours}

On an operational level, the key planning problem of a sales representative is to select the most promising customers amongst all the customer groups summarized in Figure \ref{fig:CustLife} for a visit within the next couple of days, to assign customers to days, and to determine the sequence of customers for each day in a way that the daily contractual working time is respected.
This selection should, at the same time, aim for short-term sales and maintain long-term customer relationships.
Furthermore, mandatory visits, e.g., at customers with agreed dates, may have to be included in the tour plan.
An administrative or private appointment of the sales representative, e.g., a visit at the car repair shop or at a doctor, can be also be modeled as a mandatory visit.

The goal when planning the concrete tour is to maximize the profit, i.e., to maximize sales minus costs.
To this end, customers and prospects are selected based on a score that approximates their future sales.
As described in the next section, this score can, for example, represent the expected sales or the rank of a customer or prospect compared to the others.
However, the selection of a customer does not only depend on the score but also the geographical position:
It might be beneficial to visit a customer with a low score if it is close to a high score customer or if it can be reached ``on the way'' to a high score customer or group of customers.

Personnel costs account for the largest share of the costs of a field sales force. 
Among these costs, the fixed personnel costs, such as fixed monthly salary components, are not influenced by short term decisions.
The most important variable costs are commissions, which are usually determined as a share of sales.
For this reason, this work assumes that in the short term, maximizing the scores corresponds with a maximization of the expected profit.
In practice, travel distance as the main driver of variable vehicle costs is often minimized.
This is, therefore, considered as a secondary goal in this work.


%
%
%

\subsection{Measuring and predicting a customer's value}

From a sales analytics point of view, the key challenge in planning profitable tours is to predict which customers are promising for a visit within the next days.
To this end, different data sources are used.
Promising prospects are identified based on general data often provided by the marketing department, e.g., the size of a shop, the selling potential or demographic development of the area in which the shop is located.
Prospects who are similar to profitable customers of the company are considered as promising.

For assessing the value of existing customers, data of customer relationship management systems are used, e.g., the average order volume, the time since the last order, the frequency of orders, the last visit, or historic information about interactions or information if a sales opportunity was won or lost.
The way how the value of a customer is estimated or predicted strongly varies and depends on the available database.



The effort and data that is necessary for assessing a customer's value depend on (1) the level of aggregation and (2) the content of information of the score.
The more disaggregated the level is, the more effort and data are necessary.
For example, ranking the value of three customer classes requires less data than ranking the customers or even the employees of customers on an individual level -- which is done for example in the pharmaceutical industry in the United States \citep{Mantrala2010}.
Similarly, the content of information of a score influences the effort and data needed. 
A ranking of customers based on their historical value, for example, requires less data than the prediction of sales as a response to a concrete offer for a new product.




%
%
%
%

\section{Related work}\label{sec:related_work}

In field sales force management a common procedure is to assign sales territories to sales representatives for establishing long-term customer relationships. 
The corresponding territory design problems has received significant attention in the last decades (see, e.g., \cite{Ronen1983, Rios2013, Kalcsics2019}).
The short-term profitable tour planning problem found less consideration in literature.
It relates to two lines of research:
First, we give a short overview of tour planning approaches strongly related to the problem at hand.
We continue by presenting recent approaches from literature dedicated to approximate the value of customers and prospects or to predict their future sales in the B2B context.
We finish this literature review with approaches considering both tasks.

\subsection{Tour planning}\label{sec:orienteering_literature}

The operative planning problem of selecting and routing a subset of customers such that an overall score measure is maximized has been considered in literature in the form of the orienteering problem (OP).
If tours for more than one person within the same area are planned in an integrated fashion, the problem is referred to as team orienteering problem (TOP).
The OP and its variants are subject of many publications from different domains.
Therefore, we limit this study to problem variants that are closest to our setting.
For a broad overview of related problems we refer to recent reviews by \cite{Vansteenwegen2011} and \cite{Gunawan2016}.

The planning problem considered in this work is closely related to multi-period extensions of the orienteering problem in the context of sales representative planning as well as tourist trip planning. 
Common to these problems is that tours are planned over several days  (e.g., one week) with the objective of maximizing some profit measure, typically considering working time restrictions and the inclusion of mandatory visits. 
The model that is closest to this study is the multi-period orienteering problem with multiple time windows (MuPOPTW) proposed by \citet{Tricoire2010}. 
Existing customers are modeled as mandatory visits, whereas potential customers are optional visits with an expected profit collected when they are included in a sales representative's tour.
In contrast to our setting, customers have time windows which may differ for different days in the planning horizon.
The objective is to determine one tour for each day such that all mandatory customers are visited, total profit is maximal, and daily working time restrictions, as well as a limit on the total workload in the planning horizon, are respected.
The multi-period personalized tour design problem \citep{Kotiloglu2017}, which is dedicated to planning tourist trips, similarly considers tours planned over several days with the additional consideration of opening times.
The user can specify mandatory visits that have to be included in the tours and that are complemented by optional visits associated with profits. 

For solving orienteering problems in practical applications within acceptable computational times, several authors have proposed heuristic methods, particularly population-based approaches and metaheuristics.
Population-based approaches, e.g., the particle swarm optimization inspired algorithm (PSOiA) introduced by \citet{Dang2013} and the pareto mimic algorithm (PMA) proposed by \citet{Ke2016}, have yielded excellent results even for large-scale instances of the team orienteering problem (TOP).
However, computation times remain impractical for industrial applications.
Better trade-offs between solution quality and computation times are achieved by the multi-start local-improvement approach (MS-LS) introduced by \citet{Vidal2015} and the two-phase multi-start adaptive large neighborhood search (2MLS) by \citet{Glock2020}.



\subsection{Customer scoring}
\label{subsec:value_sales}

Despite the rich body of work on orienteering problems and the impact that the scores  have on the final solutions, authors rarely consider the origin or the nature of the applied score. 
To address this gap, we provide a comprehensive overview of scoring techniques that are applied in practice, and that were recently proposed in literature.
Based on this summary, we derive modeling variants {for the profitable sales representative tour problem} in Section \ref{sec:model} and study the impact of different types of scores in Section \ref{sec:evaluation_models}. 

\citet{Bose2009} consider the selection of target customers as the ``core activity'' of direct marketing, that means of all types of marketing activities that are personalized and directed to an individual customer.
For selecting promising customers, usually, a score on the individual customer level is generated.
The extensive review of \citet{Bose2009} deals with the literature of the early 2000s for the B2C context and shows the wide range of applied techniques covering very basic statistical techniques up to machine learning or hybrid techniques.

In the B2B context, inputs and methods to estimate the scores and the scores themselves vary even more depending on industries, the size of the customer base, and available data.
Although B2B applications found less attention in literature \citep{Jahromi2014, Martinez2020}, we restrict our overview to illustrative prediction methods very recently published in literature (after 2010) and established methods described in textbooks for practitioners.
Our main goal is to derive typical categories of scores as a basis for the models introduced in the next section.
Table~\ref{tab:scores} gives an overview of the scoring methods, summarizing the required input data, applied prediction method, the obtained output and main area of application.

\begin{sidewaystable}
\renewcommand{\arraystretch}{1.05}
  \centering
  \resizebox{\textwidth}{!}{
    \begin{tabular}{p{2em}p{10em}p{20em}p{24em}p{20em}p{6em}p{10em}}
    \toprule
          & method type & type of input variables & approach & type of score & applicable to & reference, e.g. \\
    \midrule
    (1) & ABC analysis & sales, profit contribution of customers & ranking by sales or profit and classification of customers by ABC scheme (20\% of customers produce 80\% of sales and so on) & classes of customers (top, developable, unsatisfying) & customers & \citep[pp. 43 ff.]{Behle2014} \\
    \midrule
    (2) & customer value & customer statistics, interaction, sales history, etc. & human assessment in the form of scorecards with (scores e.g. -20,..,20); sum of these scores & ranking of customers by score & customers, prospects & \citep{Duncan2015}; \citep[pp. 43 ff.]{Behle2014} \\
     &  & contribution of customer in the past (toward profits), discount rate & sum of contributions of customer in the past adapted by discount  & ranking of customers by past customer value (PCV) & customers & \citep[pp. 111 ff.]{Kumar2018} \\
    &  & recurring sales and costs, lifetime of customer, discount rate, acquisition cost, externalities & different formulas, partly including prediction e.g. of the future margin, the lifetime or frequency & ranking of customers by customer lifetime value (CLV)  & customers & \citep{Venkatesan2004}; \citep[pp. 112 ff.]{Kumar2018} \\
    \midrule
    (3) & profiling based on customer classes or segments & statistics available for customers and prospects, value for base customers or clusters of customers & select prospects most similar to a set of customers (or sets of customers), e.g. the 20\% most valuable customers of the base, using for example the k-nearest neighbor algorithm  & set(s) of prospects & prospects (customers) & \citep{Dhean2013}; \citep[pp. 118 ff.]{Kumar2018} \\
    \midrule
    (4) & predicting propensities / probabilities & customer, product, sales history & logistic regression & win propensity & customers & \citep{Yan2015} \\
          &       & customer, product, interaction sequence, sales history & two-dimensional Hawkes processes & win propensity & customers & \citep{Yan2015hawk} \\
          &       & sales history, sales team, customer, interaction, behavior & logistic regression, decision tree, random forest, XGBoost & win propensity & customers & \citep{Mortensen2019} \\
          &       & sales history and derived features & lasso regression, extreme learning machine, gradient tree boosting & win propensity & customers & \citep{Martinez2020} \\
          &       & product, sales team, customer, interaction, behavior, sales history  & naive Bayes, decision tree, support vector machine, neural network; methods to explain classification models and predictions & win propensity, charts explaining influence of attributes & customers & \citep{Bohanec2017} \\
          &       & customer, product, sales history & random forest, support vector machine, XGBoost, CatBoost & win propensity, conversion probability & customers & \citep{Eitle2019} \\
          &       & sales history & different probability formulas weighting probability observations of different periods & conversion probability & customers & \citep{Xu2017} \\
          &       & online news texts, company & different survey types to collect preferences of experts about importance of decision factors, Bayes network & conversion probability & prospects & \citep{Benhaddou2017} \\
          &       & recency, frequency, monetary value & decision tree variants, logistic regression, boosting & churn probability & prospects & \citep{Jahromi2014} \\
          \midrule
    (5) & predicting win propensity and expected sales & customer, sales history & logistic regression for win propensity and mixture of expert supported workshop, quantile modelling, linear regression, and k-nearest neighbor to estimate expected sales & win propensity, expected sales & customers, prospects & \citep{Lawrence2010, Lawrence2007} \\
        & & customer, contact, behavioral inputs & boosted tree classifier for probabilities & win propensity, conversion probability, expected sales afterwards converted to rates (A, B, C, D) & customers & \citep{Duncan2015} \\
    \midrule
    (6) & customer response function & client segments, deal size, sales team, sales history & linear regression & expected sales of customer segments depending on effort by sales role & customers & \citep{Kawas2013} \\
          &       & sales history & two models established in reliability engineering using regression techniques for parameter estimation & retention probability of customer type depending on number of visits by sales representatives & customers & \citep{Golalikhani2013} \\
          &       & product, customer, interaction, sales team, sales history & multilayer perceptron neural network & win propensity for product depending on sales team,  expected sales for product & customers & \citep{Bischhoffshausen2014}; \citep{Bischhoffshausen2015}; \\
          \bottomrule
    \end{tabular}}
  \caption{Selection of typical B2B customer scoring approaches from textbooks for practitioners and published after 2010.}
    \label{tab:scores}
\end{sidewaystable}

\paragraph{Input variables}
In the following, we summarize types of input variables that are used amongst others by the scoring approaches referred to in Table~\ref{tab:scores}:

\begin{itemize}
	\item Customer or prospect profile: static information about customers or prospects, e.g., location of company, location of site, industry, sector, number of employees, revenue, market value, number of job openings, funding
	\item Contact person(s): information about the direct contact person(s) of the customer, e.g., name, title, authority level (e.g. low, medium, high)
	\item Product: information about the product, the campaign or the service, e.g., product type, (estimated) deal size, complexity of service, type of campaign 
	\item Sales representative or team: information about the sales representative or the team of sales representatives, e.g., name, skill level, role (e.g. technical, sales)
	\item Customer behavior: information about the behavior of the customer, e.g., usage of web services, response to marketing actions, information requested, proposal requested,  expression of interest
	\item Interaction: information about the interaction between companies, e.g., visit of the sales representative, phone call, information sent, product or service demonstrated
	\item Sales history: information about the sales history of a customer, e.g., history of transactions, history of sales status, average inter purchase time, time since the last order (recency), frequency of orders in a certain period, average amount per order (monetary value), cost for sales and marketing, won and lost sales opportunities
\end{itemize}

Typical sources for these data are internal systems such as customer relationship management, enterprise resource planning, or web tracking systems for the own web sites and services. Additional sources are databases of marketing agencies, e.g., about demographics and company profiles, or analytics tools for news or social media.

\paragraph{Scores and methods} 
The content of information and the level of detail of the scores described on the rows of Table~\ref{tab:scores} vary strongly:
The result of the ABC analysis in row (1) consists of three classes of customers (e.g., top, developable, unsatisfying) based on their contribution to the sales volume or profit in the past.
This approach follows the reasoning that scarce resources should be concentrated on customers who are the most profitable.
The score of an individual customer simply corresponds to the respective class A, B, or C he or she falls into.

The methods subsumed under ``customer value'' in row (2) try to estimate the value of customers for the company as a basis for the selection.
In the first method, sales representatives assess each customer or prospect based on a scorecard.
That means, they assign a score, e.g., a value between -20 and 20, to each customer and category. 
The total score is built by the (weighted) sum of scores.
Categories and scores are ``hand-tuned by experienced members of the marketing or sales team'' \citep{Duncan2015} and comprise monetary and qualitative aspects as well as expected developments of the customer.
Depending on the categories the scorecard method can also be applied to prospects.
The customer past value (CPV) and customer lifetime value (CLV) both rely on available sales and cost data for each customer:
The CPV represents the value a customer had in the past in terms of the contribution to the profit and assumes that this value is a surrogate measure for the future value.
For determining the CLV, there exist different formulas.
Usually, they additionally consider the prediction of the behavior of customers, such as the order frequency or lifetime \citep[see e.g.][pp. 112 ff.]{Kumar2018}. 
CLVs are calculated based on these independently estimated values as well as on past sales and cost data.

The profiling methods referred to in row (3) follow the idea, that prospects similar to profitable or high potential customers are more likely to become profitable customers.
Therefore, prospects are assigned to the customer class they are most similar to, or they are ranked by similarity to the most profitable customer class.
There exist different methods to segment or cluster customers \citep[see e.g.][]{Bose2009, Kumar2018}.
A possible method to assign prospects to similar classes is the k-nearest neighbor algorithm \citep[e.g][p. 354]{Linoff2011}.
The typical result of a profiling step is a group of prospects considered as profitable. 
Additionally, the similarity to a profitable customer class can be used as ranking criterion \citep{Dhean2013}. 

The methods of row (4) are dedicated to predicting probabilities that some event occurs.
To keep the acquisition process of many different prospects and customers manageable, sales teams often organize the status of customers in so-called sales pipelines or sales funnels.
A sales pipeline is a way to describe the process by dividing it into different stages \citep{Dhean2013}. 
The stages and their names are often customized by sales teams and terms have different meanings.
To give an example from our study, in \citep{Dhean2013} the stages ``suspects, prospects, leads, customers'' are used, in \citep{Duncan2015} the stages are ``awareness, leads, marketing-qualified leads, sales-qualified leads, won''.
However, the underlying logic is the same: 
If a customer or prospect reaches the next stage because the contact person responded to some marketing activities such as an e-mail, the probability to win a deal increases.
Hence, all methods subsumed in row (4) predict the probability of a customer to reach the next stage or directly the last stage.
The latter is also referred to as win propensity, which is interpreted as the probability to win a deal or a customer. 
Conversely, it is also of interest to predict if a customer has the tendency to change the supplier (customer churn).
The listed methods use all types of inputs in different combinations, and many different classification techniques from statistics and machine learning, such as logistic regression, tree-based methods, or neural networks.

The authors of \citep{Lawrence2010, Lawrence2007} and \citep{Duncan2015} in row (5) predict both, the win propensity and the expected sales in case of a won deal.
For predicting the win propensity the authors used logistic regression and a boosted tree classifier respectively.
\citet{Lawrence2010, Lawrence2007} describe the prediction of sales as a result of expert workshops followed by several quantitative methods, while \citet{Duncan2015} do not mention the method they applied.
In the case of \citet{Duncan2015}, the resulting win propensity and expected sales values are not communicated to the sales representatives. 
Instead, rates (A, B, C, D) calculated based on these values are communicated.

Customer response functions as referred to in row (6) represent how ``sales at relevant levels of aggregation vary with selling activities'' \citep{Mantrala2010}.
That means, they try to use historic data to predict the behavior of a customer in response to some concrete effort e.g., depending on the assignment of a concrete team of sales representatives to a customer \citep{Bischhoffshausen2014, Bischhoffshausen2015}, depending on the time a sales role spends to an individual customer \citep{Kawas2013} or depending on the number of visits at a customer type \citep{Golalikhani2013}.
The methods applied to estimate these response functions reach from basic statistics to estimate independent variables of a given response function, over linear regression, to neural networks.

In general, all types of scores introduced in Table~\ref{tab:scores} are applicable for selecting customers for operative tours.
They differ considerably in the content of information and the level of aggregation.
Customer classes from an ABC analysis are very coarse scores but easy to generate.
Customer value scores also reflect the future value of a customer depending on the input they are based on.
However, they are not necessarily correlated to the short term value of a visit to a customer which might be dedicated to informing about a specific product or a campaign.
In contrast, conversion probabilities or win propensity scores usually consider a concrete product or campaign and are comparable over several stages of a sales pipeline.
However, in case of strongly deviating sales volumes per customer, it might be necessary to additionally consider the expected sales volume.

Customer response functions additionally consider the effort assigned to a customer and can be considered as the most detailed method for scoring.
For example, the score functions of \cite{Bischhoffshausen2014, Bischhoffshausen2015} represent the probability to win an opportunity depending on the opportunity itself, the assigned sales team, and the customer. 
However, these methods require very detailed data of high quality and a long history.
According to \cite{Bischhoffshausen2016}, this can only be achieved if a ``rigorous data sales opportunity assessment process'' is established, and a special team only works to maintain data quality.

%
%
%

\subsection{Combination of sales analytics and planning}

We close this section with an overview of approaches that explicitly address the combination of sales analytics and optimization.
They use objective functions that are similar to profit maximization but solve different optimization problems such as territory design without explicitly considering tour planning aspects.

The joint problem of assigning customers to sales representatives' territories and of allocating visit time to customers such that expected profits are maximal is studied by \citet{Skiera1998}.
The model integrates sales response functions, but, due to a lack of accurate data, the authors use ``subjective judgments by the management'' to define the model used in their case study.  
The authors show that their strategy outperforms approaches targeted towards balanced territories even when the sales response function includes relatively large errors.
The authors do not examine the impact of these errors on the quality of the resulting solutions.
\citet{Kawas2013} address the tactical problem of assigning sales representatives to sales opportunities such that expected profits are maximal, based on empirical sales response functions.
To our knowledge, this is the first publication discussing in detail how inaccuracies in the analytics component impact the outcome of the optimization approach.

\citet{Golalikhani2013} describe an assignment problem considering spatial aspects.
They distinguish between two activities in sales force planning: Visits at established customers and visits at prospects. 
{Their objective in a first stage is to determine the number of visits to existing customers that maximize expected profits.
These profits are modeled in form of a response function that assumes that the probability for retaining a customer increases with the number of visits.
Based on the results of this phase, they determine the minimal number of sales representatives in each territory and the assignment of activities to individual sales representatives.}
While the authors discuss structural changes in the solutions due to different input parameters, they do not address the robustness towards input errors, assuming that the available information is accurate.

Finally, \citet{Bischhoffshausen2015} address the problem of assigning teams of technical and sales specialists to customers.
The probability of a sale depends on the assigned sales team and the expected sales per product is estimated.
The planning problem consists of maximizing the expected revenue, depending on the assignment of teams to customers. 
Geographical aspects are considered in a very coarse way by limiting the aggregated travel distances to customers for each sales representative.

\section{Modeling the profitable sales representative tour problem}
\label{sec:model}

In this section, we introduce models that can integrate the results of different customer value approximation methods in operative tour planning.
We first {categorize} different types of scoring methods, {based on the literature summarized in Table \ref{tab:scores}}.
We then integrate these categories in a common problem formulation.
To this end, all variants of the profitable sales representative tour problem are modeled as a multi-period orienteering problem (MP-OP).

\subsection{Profitability approximation models}\label{sec:model_derivation}

\begin{figure}[h!]
 \includegraphics[width=0.9\textwidth]{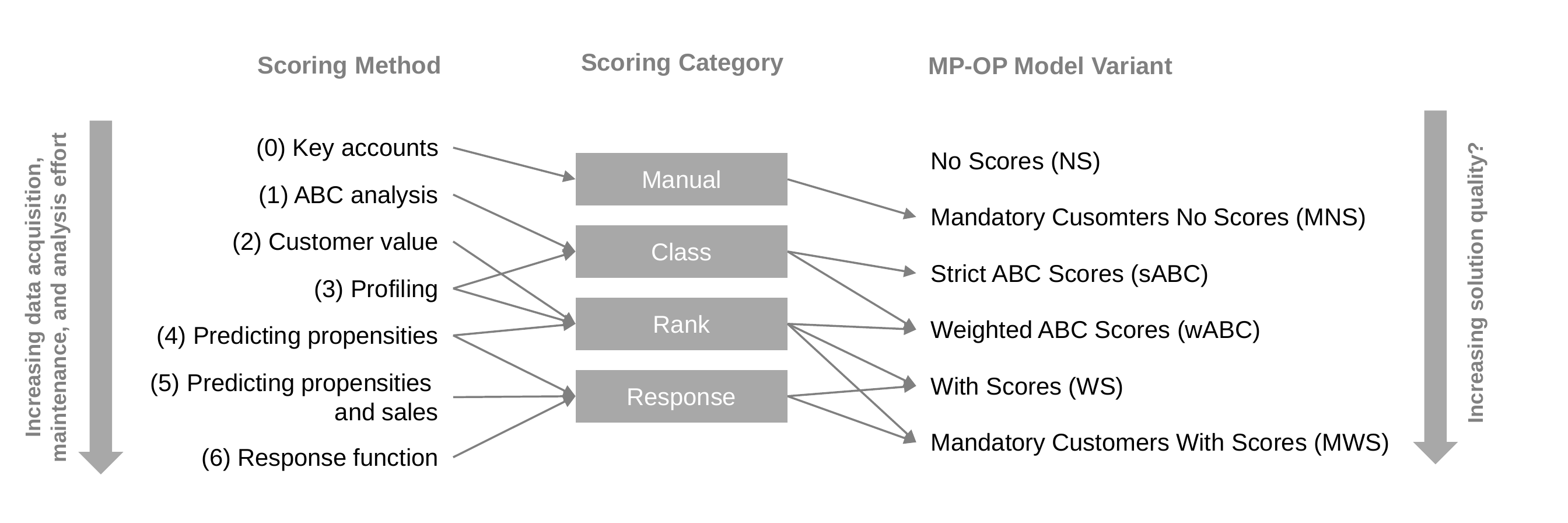} 
  \centering
  \caption{Relationship between scoring method, scoring category, and MP-OP modeling variant.}
  \label{fig:scareCat}
\end{figure}

As described in Section~\ref{subsec:value_sales}, there exist different types of scores for estimating the value of a customer or for approximating future sales of a customer.
They differ in the level of aggregation and the level of content.
{Based on the literature review and the experience from practice, we propose four categories of customer scoring methods based on their output, and introduce model variants of the MP-OP for each of them.
An overview of the relationship between the scoring methods (Section~\ref{subsec:value_sales}), the scoring categories (this section), and the MP-OP modeling variants (Section~\ref{sec:profitability_models}) is depicted in Figure~\ref{fig:scareCat}}.

{The four scoring categories can be characterized as follows:}

\begin{itemize}
	\item Manual: Even if no data-driven scoring phase is in place, sales representatives usually manually select a set of especially promising customers {(often referred to as key account customers)}, who they definitely want to visit.
	These customers build the core structure of tours if the sales representatives plan their tours without IT support.
	\item Class: All methods resulting in customer classes, such as the ones of line (1) in Table~\ref{tab:scores} but also the approach of \citet{Duncan2015} from line (5), which assigns rates (A, B, C, D) to customers, are subsumed under the type ``class''.
	\item Rank: Scoring methods are considered as ``rank'' methods if {scores assigned to customers do not necessarily have a linear relationship to the expected profit of a visit. Instead, they can be understood as an ordinal scale for the purpose of comparing customers.} 
	In this case, one can state that a customer with a score twice as high as the score of another customer should be given priority but the expected outcome of the visit is not worth twice as much.
	This applies to the methods of line (2) of Table~\ref{tab:scores}, where approximations for the value of a customer based on historical values are proposed. 
	Methods of line (4), which determine win propensities or conversion probabilities not related to a concrete product or sales volume, also belong to the type ``rank''. 
	It still might be the case that a customer with a high win propensity but a low volume would be less promising than a customer with lower win propensity but a high volume. 
	\item Response: A method is considered as a sales response function method if it estimates the direct response of a customer to a concrete visit presenting the current offers.
	 That includes, for example, a prediction of the estimated (additional) sales.
	 We talk about additional sales because in some industries, customers purchase regularly and sales representatives try to generate additional sales by presenting the newest offers. 
	 Hence, the methods of the last line of Table~\ref{tab:scores} are subsumed under ``response'' scoring methods.
	 However, depending on the data input also methods from line (4) and (5) might be used as response function method:
	 The expected win propensity for a concrete product (e.g., a new license) with a fixed price relates to the expected profit in a linear way.
\end{itemize}

In the following, we propose model variants that allow the use of these four categories of scores in operative tour planning. 
We do so by transferring these categories into different parameterizations of the objective function of the MP-OP.

\begin{center}
\begin{table}[htbp]
\begin{tabularx}{\textwidth}{l X} 
\toprule
Sets and parameters & \\
\midrule
$\mathcal{N}^C$ & Set of customer nodes \\
$\mathcal{N}$ & Set of nodes including customer nodes $\mathcal{N}^C$ and home location (node $0$) \\
$\mathcal{N}^{M}$ & Set of mandatory customers $\mathcal{N}^{M} \subseteq \mathcal{N}^C$ \\
$D$ & Set of considered working days\\
$t_{ij}$ & Time for traveling from node $i \in \mathcal{N}$ to node $j \in \mathcal{N}$ \\
$s_i$ & Service time associated with one visit at node $i \in \mathcal{N}^C$  \\
$p_i$ & Score of customer node $i \in \mathcal{N}^C$ \\ 
$T$ & Maximum working time per day \\
\midrule
Decision variables & \\
\midrule
$x_{ijd} \in \{0, 1\}$ & Equals 1 iff node $i \in \mathcal{N}$ is visited directly before node $j \in \mathcal{N}$ on day $d \in \mathcal{D}$\\
$v_{id} \in \{0, 1\}$ & Equals 1 iff node $i \in \mathcal{N}^C$ is visited on day $d \in \mathcal{D}$\\
\bottomrule
\end{tabularx}
\caption{Sets, parameters and decision variables for the MP-OP.} 
\label{table:mp-op}
\end{table}
\end{center}

\subsection{The multi-period orienteering problem}

The common base model for planning profitable tours for all customer scoring variants is the MP-OP.
Based on the sets, parameters, and decision variables given in Table~\ref{table:mp-op}, this model can be formulated as follows:
\begin{align}
	\sum_{i \in \mathcal{N}^C \setminus \mathcal{N}^M}\sum_{d \in \mathcal{D}} v_{id}~p_i \quad \rightarrow \quad \max \label{e:bm0}
\end{align}
s.t.
\begin{align}
&\sum_{d \in \mathcal{D}} v_{id}  \leq 1 & i \in \mathcal{N}^C  \label{e:bm1} \\
&\sum_{j \in \mathcal{N}} x_{ijd} = \sum_{j \in \mathcal{N}} x_{jid} & i \in \mathcal{N}, d \in \mathcal{D} \label{e:bm2}\\
&  \sum_{j \in \mathcal{N}} x_{0jd} = \sum_{j \in \mathcal{N}} x_{j0d} = 1 & d \in \mathcal{D} \label{e:bm3}\\
& \sum_{j \in \mathcal{N}} x_{ijd} = v_{id} & i \in \mathcal{N}, d \in \mathcal{D} \label{e:bm4}\\
& \sum_{i,j \in \mathcal{N}} (t_{ij} + s_{i}) ~ x_{ijd} \leq T & d \in \mathcal{D} \label{e:bm5}\\
& \sum_{d \in \mathcal{D}} v_{id} = 1 & i \in \mathcal{N}^{M} \label{e:bm6}\\
& \sum_{i,j \in \mathcal{Q}} x_{ijd} \leq |\mathcal{Q}| - 1 & \mathcal{Q} \subseteq \mathcal{N}^C, d \in \mathcal{D} \label{e:bm7}\\
& x_{ijd} \in \{0,1\} & i,j \in \mathcal{N}~|~i \neq j, d \in \mathcal{D} \label{e:bm8}\\
& v_{id} \in \{0,1\} & i \in \mathcal{N}^C, d \in \mathcal{D} \label{e:bm9}
\end{align}


The objective function in equation (\ref{e:bm0}) maximizes the sum of realized scores, i.e. the sum of scores of the customers, which are selected.
Equations (\ref{e:bm1}) ensure that each customer $i$ is visited at most once.
The next set of constraints ensures flow conservation, while constraints (\ref{e:bm3}) ensure that the tour starts and ends at the home location of the sales representative.
With equations (\ref{e:bm4}) we link the decision variables $x$ and $v$, i.e. we ensure that, if a customer $i$ is visited on day $d$, an outgoing arc is selected for that customer on that day and the other way around.
With constraint set (\ref{e:bm5}) the maximum working time per day $T$ is enforced, considering both travel time and service time.
Equations (\ref{e:bm6}) make sure that all mandatory customers are visited once within the planning horizon.
Constraint set (\ref{e:bm7}) eliminates subtours.
Sets (\ref{e:bm8}) and (\ref{e:bm9}) restrict the decision variables to be binary.

\subsection{Modelling variants for representing customer profitability estimates}\label{sec:profitability_models}

This section formally defines the model variants that allow us to represent the different categories of profitability approximation methods summarized in Section \ref{sec:model_derivation}. 
Table~\ref{tab:mp-op} gives an overview how the results of different scoring categories are translated into objective functions for the tour planning phase.
The modeling variants based on these objective functions are introduced and further discussed in the remainder of this section.

\begin{center}
\begin{table}[htbp]
\begin{tabularx}{\textwidth}{l l X} 
\toprule
Scoring category & Model variant & Description \\
\midrule
No & NS & Every customer is treated equally. \\[0.5em]
Manual & MNS & Manually selected customers are treated as mandatory customers. Remaining customers are treated equally. \\[1.5em]
Class & sABC & Customers in a higher class are strictly preferred to customers of a lower class. \\[1.5em]
Class / Rank & wABC & Customers of the same class are evaluated by the average score of the class or a given number of customers in a lower class substitute the visit of a customer of a higher class. \\[1.5em]
Rank / Response & WS & Each customer is evaluated based on the respective score. \\[1.5em]
Rank / Response & MWS & Same as WS. Mandatory customers are added, e.g. the ones with the highest scores or the ones who declared an intent to purchase in a personal contact. \\
\bottomrule
\end{tabularx}
\caption{Overview of scoring category and MP-OP variant.} 
\label{tab:mp-op}
\end{table}
\end{center}

\paragraph{No Scores (NS)}

In the most simple version, no scores for customers are considered.
If no scoring phase is in place or no data about customers are available, each customer has the same value and the MP-OP maximizes the number of visited customers over the planning horizon.
In this case, we set the score parameter to $ p_i = 1$ for all $i \in \mathcal{N}^C$.
Since the contribution to the objective function is equal for each customer, the selection of a customer only depends on the geographic position in relation to the other selected customers.

\paragraph{Mandatory Customers No Scores (MNS)}

The MNS model relies on the manual selection of the most promising customers:
These customers are added to the set of mandatory customers $\mathcal{N}^{M}$, and constraint set (\ref{e:bm6}) ensures that they are part of the resulting tour plan.
As in the NS model, each customer has a score of $p_i = 1$ for all $i \in \mathcal{N}^C$. 
This means, that besides all mandatory customers as many remaining customers as possible are selected for a visit depending on their geographical position.

\paragraph{Strict ABC-scores (sABC)}

If customer classes are the result of the scoring phase, the score parameters $p_i$ are calculated based on the class of each customer $i \in \mathcal{N}^C$.
There exist two ways to calculate these parameters: strict and weighted ABC-scores.

Strict ABC-scores represent a hierarchical objective function:
If it is possible, all A customers are visited before a B customer is added to a tour, and all B customers are visited (if possible) before a C customer is selected.

If the customer classes are given as sets of customers $\mathcal{N}^{AC}, \mathcal{N}^{BC}, \mathcal{N}^{CC}$ (with $\mathcal{N}^{AC} \cup \mathcal{N}^{BC} \cup \mathcal{N}^{CC} = \mathcal{N}^{C} $) representing the set of A, B, and C customers respectively, this can be realized based on the following score parameters $p_i$:
\begin{align}
	p_i = p^C  & = 1 & i \in \mathcal{N}^{CC} \\
	p_i = p^B  & = p^C \cdot |\mathcal{N}^{CC}| +1 & i \in \mathcal{N}^{BC} \\
	p_i = p^A  & = p^C \cdot |\mathcal{N}^{CC}| + p^B \cdot |\mathcal{N}^{BC}| +1 & i \in \mathcal{N}^{AC}
\end{align}

\paragraph{Weighted ABC-scores (wABC)}

The wABC model is applicable for scoring methods of type ``class'', if an average class score can be determined or if the relative importance of customers in a class is assessed by the sales representative.

If the classification of customers is the result of a two-step approach, in which a score is calculated for each customer, and afterward the customers are classified based upon this score as reported e.g. by \citet{Duncan2015}, the score value $p_i$ of customers of the same class can be set to the average score of the class.

If the classification is not the result of a two-stage process, the  value of customers within one class must be based on a judgmental assessment:
The sales representative estimates how many B customers can substitute the visit of an A customer and so on.
If, for example, the visit of three B customers is considered as valuable as the visit of one A customer, and five C customers substitute one B customer, the scores are set as follows:
\begin{align}
	p_i = p^C & = 1 & i \in \mathcal{N}^{CC}\\
	p_i = p^B & = 5 \cdot p^C = 5 & i \in \mathcal{N}^{BC} \\
	p_i = p^A & = 3 \cdot p^B = 15 & i \in \mathcal{N}^{AC}
\end{align}

For the experimental study in this work, we used the values $p^C=1$, $p^B=5$, and $p^A=15$, a possible judgmental assessment of a sales representative.
Note that both models, sABC and wABC, can be easily extended to variants with more than three customer categories.

\paragraph{With Scores (WS)}

If a ``response'' method was applied in the scoring phase, each customer is evaluated in the tour planning phase with its resulting score.
In this case, the parameter $p_i$ is simply set to the predicted (additional) sales or to the predicted (additional) win propensity for each customer $i \in \mathcal{N}$.

\paragraph{Mandatory Customers With Scores (MWS)}

In some cases, both, customer scores and mandatory customers are considered.
This might be the case if scores are calculated with a response function and the sales representative has additional information, e.g., about the intention of customers to place an order.
To assure that these customers are visited, they are added to the set of mandatory customers, while parameter $p_i$ is set to the predicted (additional) sales or to the predicted (additional) win propensity of each customer $i \in \mathcal{N}$.

Note that considering mandatory customers can have the effect that instances might be infeasible.
However, in practice, the number of mandatory customers is usually rather small compared to the total number of customers.
{Nevertheless, if infeasibilities occur, mandatory customers can be treated as optional customers with a very high score resulting in a similar model behavior as the sABC formulation reduced to two classes.}

\subsection{Solution approaches}
\label{sec:solap}

Only small instances of the MP-OP with less than 20 customers can be solved to optimality within reasonable run time limits by solving the presented formulation with a state-of-the-art solver for integer linear programs (see Section~\ref{sec:ev2mls}).
{For operative tour planning components, users of PTV solutions accept run times no longer than five minutes.
Hence, a heuristic approach becomes necessary for solving the MP-OP on real-world instances with several hundreds of customers.}
We apply the two phase multi-start adaptive large neighborhood search (2MLS), a framework for solving team orienteering problems (TOP) with different objective functions first introduced in \citet{Glock2020}. 
Compared to state-of-the-art methods from the literature, 2MLS shows competitive results with a very good trade off between run time and solution quality on a broad range of TOP benchmark instances \citep{Glock2020}.
{For the interested reader,} we briefly describe the main components of the 2MLS {and how we tailored them to the problem at hand} in Section~\ref{sec:2mls} of the Appendix.

\section{Test instances from industry}\label{sec:instances}

For evaluating different model variants and the 2MLS heuristic we used both, two real-world instance sets from clients of PTV Group and an instance set from literature.
{As none of the clients of PTV Group was able to provide us with raw data for an end-to-end evaluation, including the prediction stage, and, to the best of our knowledge, there exists no corresponding data set in literature, we adapted and supplemented the original instances.}
We furthermore introduced new sets based on the base instances by varying some important characteristics for the computational studies with a special focus on prediction and classification errors.
{These synthetic data sets setup allows us to evaluate the influence of the customer approximation models on the results, including the impact of potential prediction errors.}

In the following, we introduce the original instance sets, their characteristics, as well as the way we adapted and supplemented the data.

\subsection{Base instance sets}

In the following, our three base instance sets, \Nesp, \Brill, and \Tric, are introduced from a general perspective, with a closer look at the characteristics of the score, and the geographical distribution.
Tables~\ref{tab:setMetrics} and \ref{tab:setScores} summarize the characteristics of the three sets.

\begin{table}[htbp]
  \centering
    \begin{tabularx}{\textwidth}{p{0.06\textwidth}p{0.11\textwidth}p{0.09\textwidth}p{0.12\textwidth}p{0.12\textwidth}p{0.17\textwidth}p{0.13\textwidth}}
    \toprule
    Set & Num.\newline{instances} & Horizon\newline[days] & Num.\newline customers & Mand. \newline customers & Service time\newline [min] & Tour\newline dur. [h] \\
    \midrule
    \Nesp & 17  & 5 & 81 to 283 & 15 & 15 to 45 & 8 \\
    \Brill & 16 & 4 & 25 to 55 & 15  & 60  & 8   \\
    \Tric & 60 & 3 & 16 to 152 & 15  & 5 to 60 & 8   \\
    \bottomrule
    \end{tabularx}%
    \caption{Metrics of the base instance sets.}
  \label{tab:setMetrics}%
\end{table}%

\begin{table}[htbp]
  \centering
    \begin{tabularx}{\textwidth}{p{0.12\textwidth}p{0.12\textwidth}p{0.12\textwidth}p{0.5\textwidth}}
    \toprule
    Set      & Min. score & Max. score & Score calculation \\
    \midrule
    \Nesp      & 60 & 2.000 & Based on estimated customer ``importance'' (e.g., due to its location) and the number of days the customer is overdue with respect to a planned date \\
    \Brill     & 1 & 1.300 & No information \\
        \Tric      & 1 & 24.000 & Score corresponds to profit estimations of the company; in the study of \cite{Tricoire2010} the profits were combined with a synthetically generated probability to gain a customer\\
    \bottomrule
    \end{tabularx}%
    \caption{Scores of the base instance sets.}
  \label{tab:setScores}%
\end{table}%

\subsubsection{Application and instance characteristics}

\paragraph{\Nesp} 
The first set of instances comprises the planning problems of 17 salespersons selling small household appliances to retail shops in different areas of France.   
As shown in Table~\ref{tab:setMetrics}, the planning horizon is one week with five working days, while the maximum contracted daily working time and, hence, the maximum working time is 8 hours. 
Each sales representative serves between 80 and 280 customers over several weeks.
For each customer, the service times are estimated varying between 15 to 45 minutes.
Furthermore, for each customer a score is determined in a manual process:
It considers the customer's importance, as well as the number of days a visit is considered to be overdue with respect to some target date.
The resulting scores vary from approximately 60 to 2000.
Based on the customers' locations (street, city), the geographical location, the road distances and travel times were determined using the geocoding and route planning functions of PTV xServer\footnote{https://www.ptvgroup.com/en/solutions/products/ptv-xserver/developer-zone/ (23.01.2020)}.
Since the home locations of the sales representatives were not included in our data set, we placed them at the unweighted center of their corresponding customers.

\paragraph{\Brill} The second set of instances is from a client of PTV Group selling products to specialized building supplies stores in Germany.
The company employs 16 sales representatives, each of them serving between 25 to 55 customers.
Service times are approximately the same for all customers and are estimated to approximately one hour.
The score for each customer is given by our project partner without further explanation.
Since we have no data about the contracted working time per sales representative, and the number of customers per sales representative is rather small, for our study, we limit the planning horizon to four days with a maximum daily working time of 8 hours.
In case of a longer planning horizon, all customers can be visited with eight hours of working time and no selection of customers is necessary.
Again, based on the customer and home locations of the sales representatives, we determined geographical locations, road distances and travel times using the corresponding PTV xServer functions.
In general, the travel times are shorter than in our first set of real-world instances \Nesp.

\paragraph{\Tric}
In \cite{Tricoire2010}, the MuPOPTW along with a set of benchmark instances\footnote{{https://prolog.univie.ac.at/research/OP/}} was introduced.
The 60 benchmark instances are created by \citet{Tricoire2010} based on information of an industrial partner with five sales representatives selling software for geomarketing applications.
The original information included a set of customer locations, a travel time matrix, and potential profits for each customer.
In the resulting instances, a sales representative serves 10 to 100 optional customers and 6 to 12 mandatory customers.
In this case, mandatory customers are the ones who have to be visited regularly (e.g., every four weeks), while optional customers are prospects.
The planning horizon of this set is three days, while the service time of the customers lies between 5 to 60 minutes, and the maximum daily working time is 10 hours.
The original scores represent estimations for the profits of the customers.
For the computational study, the authors of \cite{Tricoire2010} synthetically generated a probability to gain a prospect (win propensity) and multiplied this probability with the profits to generate the score.
The instances also contain time windows for some of the customers and a total working time restriction over the three days of 24 hours.

For adapting the instances to our use case, the following changes were applied: 
We did not consider the time window and the working time restriction throughout the planning horizon.
In return, we reduced the working time limit to 8 hours. 
As score we used the profit estimations originally given by the company.
While this prevents comparison with previous results in the literature, it provides an additional example of customer importance estimates based on real-world data with different characteristics compared to our own data sets.



\subsubsection{Characteristics of scores}

\begin{figure}[h!]
  \includegraphics[width=\textwidth]{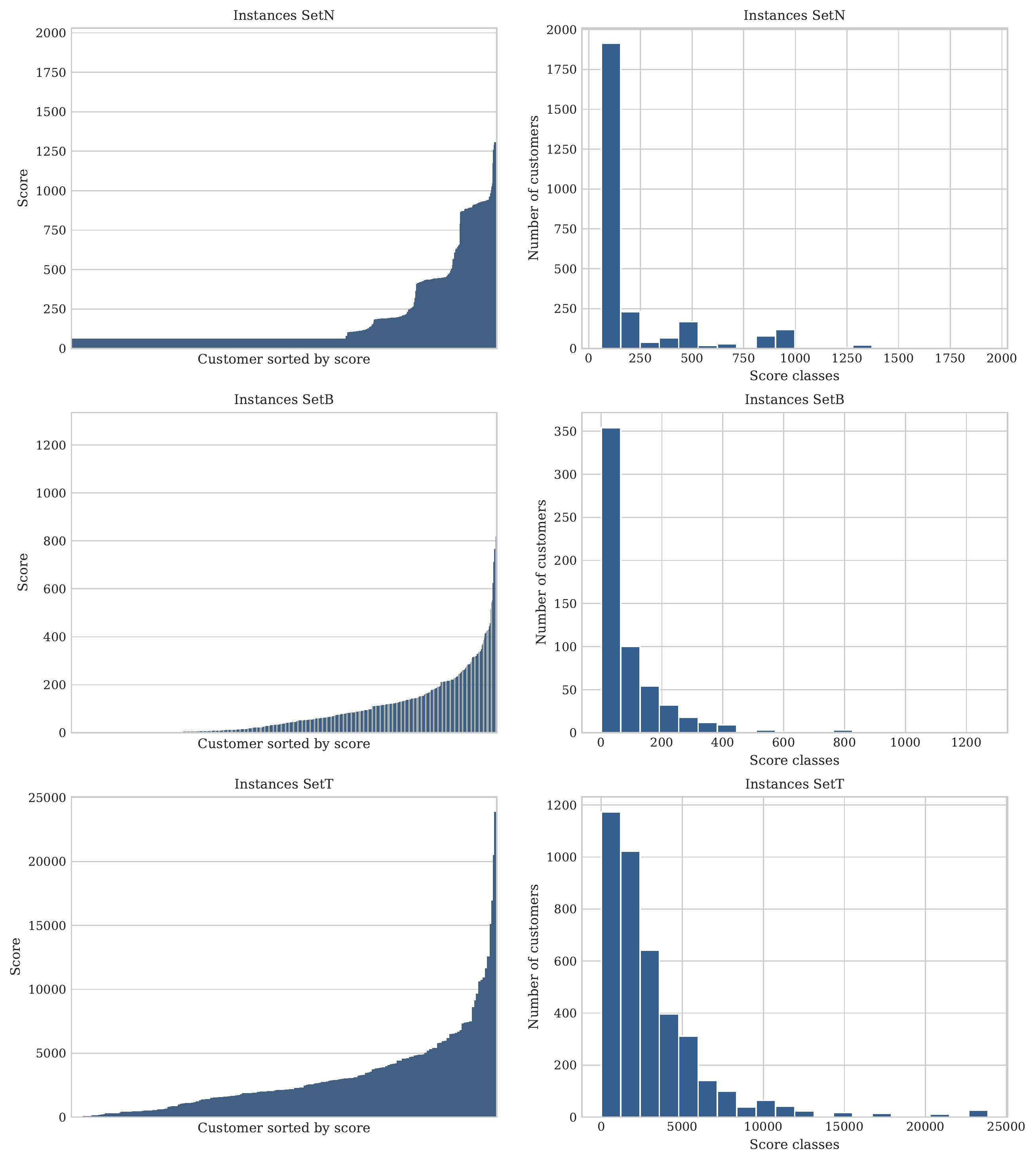} 
  \caption{Scores for all customers of all instances plotted over the sorted list of customers (left column) and histogram of scores over all customers of all instances (right column) shown for the three base sets.}
  \label{fig:PriosByBaseSet}
\end{figure}

For all three base sets the distribution of scores is shown in Figure~\ref{fig:PriosByBaseSet}. 
In the left-hand column, the score is plotted over the customers of all instances of each set, sorted in ascending order of their score value.
The right-hand column shows the distribution of scores over all customers of all instances.
The figures show that the distributions in all three cases are strongly skewed to the left, i.e. there are many customers with relatively low scores and only a few with high score values.

We do not have information about the models applied for estimating the scores, nor a sales history for the customers or even the raw data used for scoring.
In our computation study, we therefore assume that the values given are perfect (i.e. correct) forecasts for the scores. 
In a sensitivity analysis, we investigate the robustness of the optimization stage concerning errors in the prediction stage (see Section~\ref{subsec:pred_quality}).

\subsection{ABC clustering and mandatory client selection}
\label{sec:abc}
\paragraph{ABC clustering}

Since the data sets do not contain ABC classifications, we synthetically generated ABC categories for all instances of all instance sets.
To generate ABC categories consistent with the priorities for the customers of an instance, we applied a k-means clustering algorithm\footnote{We used the ``sklearn.cluster.KMeans'' implementation of the Python library scikit-learn 0.20.1 \citep{scikit-learn} with default parameter settings, i.e. with k-means++ as initialization method, with a best out of 10 multi-start runs, and a maximum of 300 iterations per single run.} using the given priority as distance measure and setting the number of resulting clusters $k$ to three.
A client who is a member of the cluster with the highest average priority is considered to be an A customer, a member of the cluster with the second-highest average priority is considered to be a B customer, and so forth.

\begin{figure}[h!]
  \includegraphics[width=\textwidth]{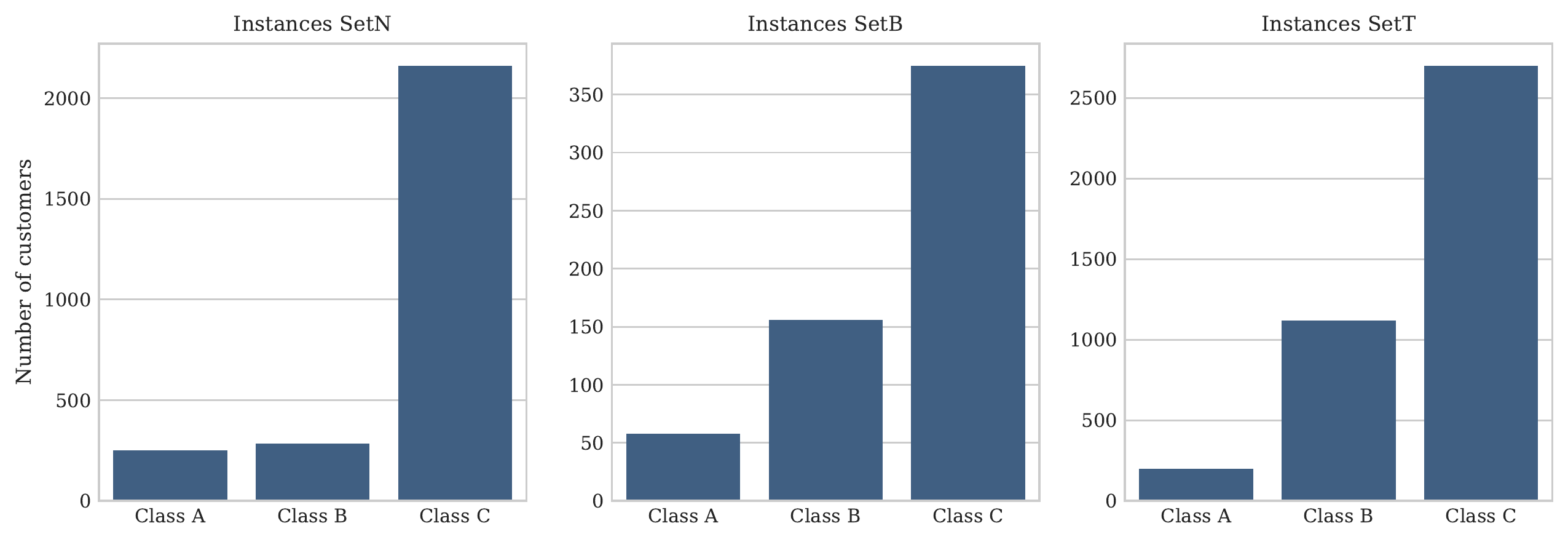} 
  \caption{Class distributions over all instances of the three base sets.}
  \label{fig:ClassBySet}
\end{figure}

As described in Section~\ref{subsec:value_sales}, a classification of customers is not necessarily based on predicted scores. 
However, for testing our models, we classify the customers based on the scores assuming that they correspond to the perfect forecast.
In practice, a wrong classification is very likely.
As shown in Figure~\ref{fig:ClassBySet}, this procedure leads to clusters differing considerably in size, if the distribution of the scores is skewed as in our three base sets. 

\paragraph{Mandatory customers}

We selected the $m=15$ customers with the highest scores as mandatory customers.
This value was chosen so that that the resulting instances are feasible in nearly all cases for all three sets.
In most cases, only A customers are selected.
However, in some instances with a low number of A customers, also B, and even C customers are selected.
Note that in practice, it might be also the case that sales representatives mark customers as mandatory irrespective of their score.
This is done, e.g., to ensure that customers considered overdue are visited, or because visits have already been confirmed to the customer and therefore have to be performed.

\subsection{Score distribution and level variation}

\begin{figure}[h!]
  \includegraphics[width=\textwidth]{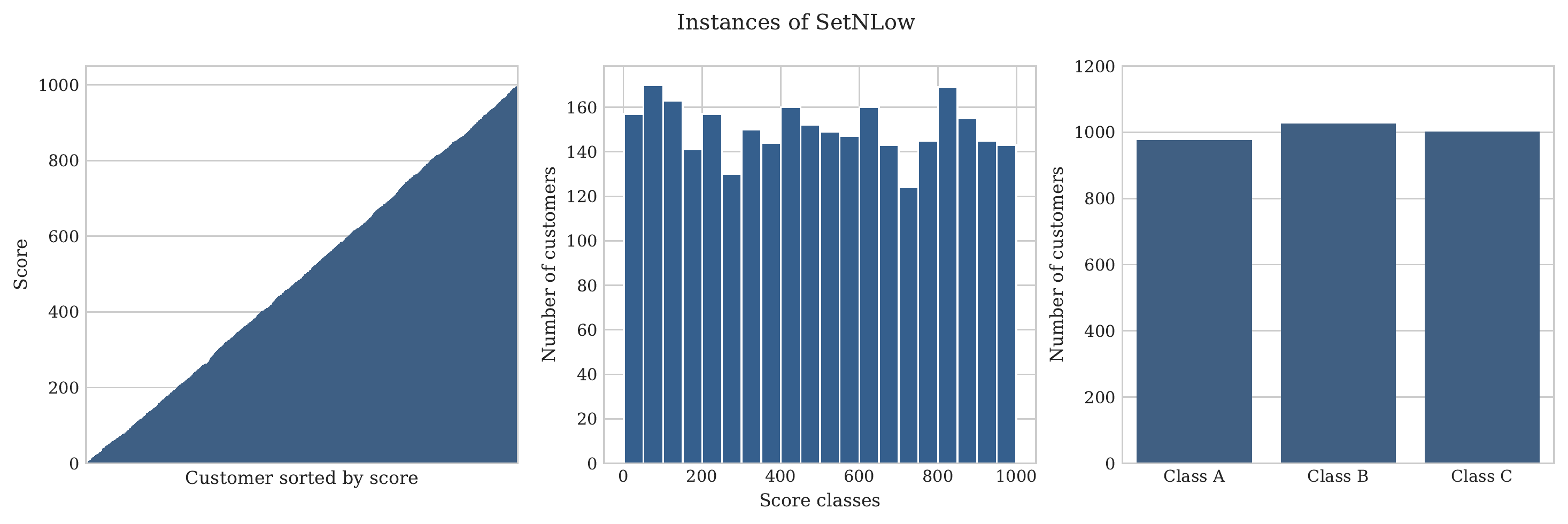} 
  \caption{Characteristics of the score distribution over all instances of \NespOne.}
  \label{fig:ScoreVariation}
\end{figure}

To study the impact of the score distribution and the differences in the level of scores, we created two more instance sets derived from each of the two sets \Nesp and \Brill:
To contrast the skewed distributions of the base instances, we generated instances assuming an equal distribution of scores between a lower and an upper bound.  
For the sets ``Low'' we assumed the lower bound to be 1 and the upper bound to be 1,000, while for the sets ``High'' we assumed a lower bound of 1 and an upper bound of 25,000.
The limits are based on the limits which appeared in the three instance sets from industry (see Table~\ref{tab:setScores}). 
The resulting characteristics of the instance set \NespOne is shown in Figure~\ref{fig:ScoreVariation}.

\subsection{Generation of small instances for exact solutions}

For solving instances to proven optimality based on the mathematical formulation described in Section~\ref{sec:model}, we generated a set of instances for which we randomly selected $n$ customers from the instance sets \Nesp and \Brill keeping the customer specific parameters of the original instances.
The selection ensured that the relative size of the sets A, B, and C  from the source instances was preserved.
The $m$ customers with the highest priority were set to be mandatory, all other customers are left optional irrespective of their status in the source instance.
To keep the capacity of the sales representative relative to the customers small enough to provoke a selection decision, we reduced the number of days per week to $d$.
We generated instances with the following two parameter sets: instances of size 10 with $\{n =10, m=2, d=2\}$ and instances of size 15 with $\{n =15, m=5, d=3\}$.


\subsection{Prediction and classification errors}
\label{subsec:pred_quality}

The predicted scores $p_i$ given from the companies are used in the computational study as input for the tour planning phase.
However, in practice, predicted scores and customer classifications contain errors.
For evaluating the robustness of our models concerning these prediction or classification errors we conduct a sensitivity analysis.
As we do not know the models nor the data used for scoring, we generate simulated score realizations $p_i^{sim}$ that represent the actual value of a customer which is revealed after a visit.
For this purpose, we assume that the models selected by the companies result in unbiased forecasts and in residuals that are independently distributed following a normal distribution with a mean value of zero and a standard deviation of $\sigma$.
Based on this assumption, we generate the simulated score realizations as follows:
\begin{align}
	p_i^{sim} = p_i + e_i \quad\quad i \in \mathcal{N}^C
\end{align}
with $e_i \sim \mathcal{N}(0,\sigma)$.
We vary the standard deviation of the errors relative to the mean score of an instance by considering different levels of the coefficient of variations $coe$.
The coefficient of variations is defined as the ratio of the standard deviation $\sigma$ and the mean score $\bar{p} = \frac{1}{|\mathcal{N}^C|}\sum_{i \in \mathcal{N}^C} p_i$.
Hence, we determine the different levels of the standard deviation by $\sigma = coe \cdot \bar{p}$.

For each instance of the instance sets SetN and SetB, we generate 10 scenarios for each of the levels $coe=[0.1, ..., 1.3]$.
A $coe$ of $0.1$ represents use cases with a very high prediction and classification quality.
Especially in case of well known customers, e.g. supermarkets ordering on a weekly basis, high quality results are possible. 
In industries such as the software sector, in which customers buy on lower frequency, and the probability that a customer positively reacts on a concrete offer are small, the prediction errors might be considerably higher.
Such a use case is represented by a $coe$ of $1.3$.

Note that we do not explicitly generate simulated classifications since the classifications are based on the predicted scores $p_i$ (see Section~\ref{sec:abc}). 

\section{Evaluating algorithmic performance of 2MLS}
\label{sec:ev2mls}

In this section, we evaluate 2MLS by comparing its results for the small instances with the ones obtained by applying a state-of-the-art integer linear programming solver to the model variants introduced in Section \ref{sec:model_derivation}.
We furthermore discuss the performance of the 2MLS on the large instances.
For generating exact solutions on the small instance sets, we modeled and solved the problem using the Python API of the mathematical solver Gurobi 8.0.1 with a MIP gap of 0.05 and a time limit of 600 seconds.
The 2MLS is implemented in C++.
{As 2MLS contains non-deterministic elements, it} is run 10 times on each instance.
The exact approach was run on an Intel core i7-8550U machine with 1.8 GHz, 16 GB RAM, and 64-Bit Windows 10, while the 2MLS was run on a Intel Xeon machine with 2.6 GHz, 119 GB RAM, and 64-Bit Windows Server 2012. 


\paragraph{Results on small instance sets}

In total we run both algorithms on 594 instances of size 10 and 594 instances of size 15. 
The exact approach solves all instances of size 10 to optimality with a computation time of approximately 1~second, or proves the infeasibility of four instances due to the mandatory customers.
For instances of size 15 the share of instances solved to optimality within the given run time limits sinks to 78\%, while 2 instances in these sets are provably infeasible.
On the 1,057 instances for which an optimal solution has been found, 2MLS finds the optimal solution for 1,042  (98.5~\%) in at least one out of ten runs, in 869 cases (82.2~\%) it finds the optimal solution in all runs.
The average gap of 2MLS to the optimal solution is 0.32~\%. 
Computation times are well below a second for all instances.

\paragraph{Results on large instance sets}\label{sec:ALNS_large}

\begin{figure}[h!]
  \includegraphics[width=\textwidth]{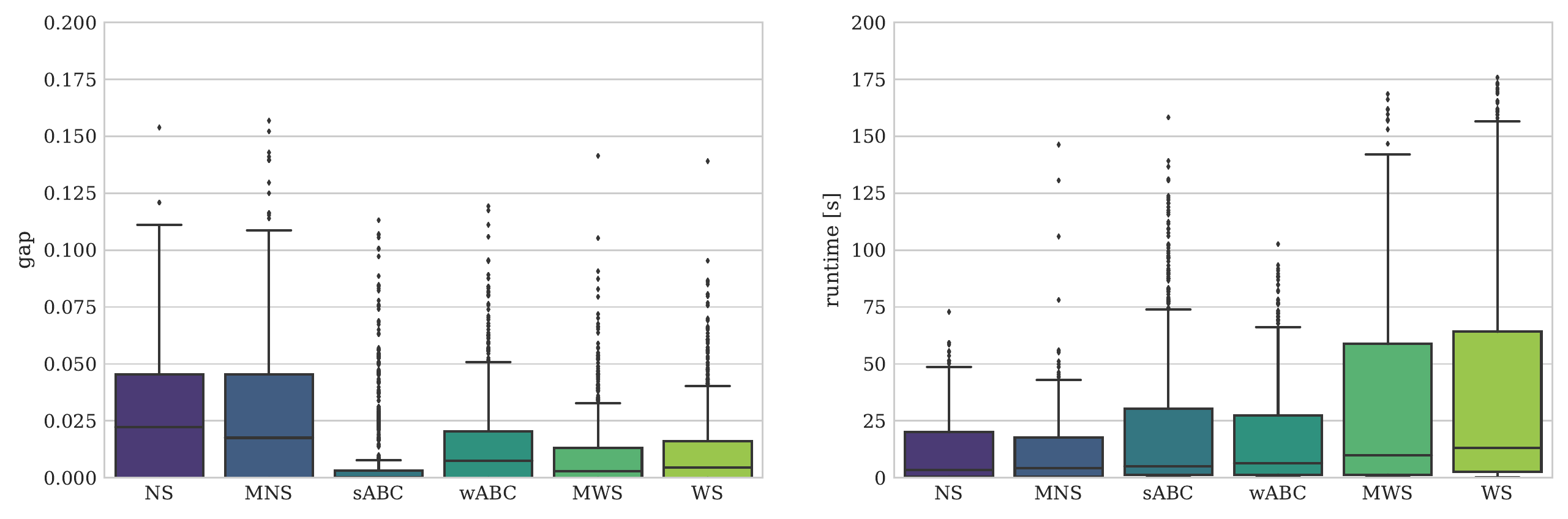} 
  \caption{Boxplots of the gap to the best found solution (left) and the runtime of solution (right) over all instance sets.}
  \label{fig:resultsALNS}
\end{figure}

Figure \ref{fig:resultsALNS} gives the deviation {between the found and the best known solution} on the left-hand side over all instances of all sets. 
The right-hand side plot gives the corresponding run times.

Even for the larger instances with up to 280 customers, run times for 2MLS are below three minutes over all model types.
Most of the instances are solved in less than a minute. 
Hence, run times of 2MLS are considered as acceptable by practitioners.

The deviation of the solution quality to the best found solution over ten runs is below 2.5~\% for a very large portion of instances in the case of the sABC, wABC, MWS, and WS models, and below 5~\% for the NS and the MNS model.
The reason is that sABC, wABC, MWS, and WS models with their dissimilarities of scores provide a structure that the search can exploit for quickly identifying the most promising solutions, and different solutions clearly dominate one another.
In contrast, NS and MNS with homogeneous scores, lead to a stronger variation of the solution quality.
In this case, solutions are highly symmetric and the score values provide no information for guiding the search.
This means that the search stops quickly due to the convergence criterion and the lack of an explicit diversification step.

In summary, we can state that 2MLS, which delivers competitive results on TOP benchmark instances (see Section~\ref{sec:solap}), also is applicable for solving profitable sales representative tour problems:
A large number of small instances is solved to optimality, run times from a few seconds up to three minutes are acceptable in practice, while the variance in solution quality over ten runs is low, especially if scores are not homogeneous.

\section{Evaluating profitability approximation models}\label{sec:evaluation_models}

In this section, we analyze the computational results with respect to the different variants for modeling the profitability of customers that were introduced in Section \ref{sec:profitability_models}.
Therefore, we mainly compare the following performance indicators of a solution:

\begin{itemize}
	\item The \textit{share of visited customers} corresponds to the number of visited customers in relation to all customers of an instance. 
	\item The \textit{share of realized score} corresponds to the realized score of a solution in relation to the aggregated score of all customers of an instance. 
\end{itemize}

{The first indicator gives an idea about whether we visit enough customers, while the second indicator gives an idea about whether we visit the profitable part of the customer base.}

To not dilute the results, we filtered the instances for the analysis as follows:
We only consider instances for which no model variant delivered a solution in which all customers can be served.
We furthermore excluded instances for which model variants containing mandatory customers (MNS, MS) become infeasible.

We start our discussion showing the impacts of the different models on optimal or near to optimal solutions for the small instance sets.
This ensures that the observed effects do not depend on the heuristic approach, but on the modeling variants that are applied.
Afterward, we analyze the solutions generated with our 2MLS approach on the large instance sets to confirm the results of the small instance sets and to show effects of the heuristic approach.

\subsection{Results on small instance sets}

\paragraph{Number of visited customers:}

\begin{figure}[h!]
  \includegraphics[width=\textwidth]{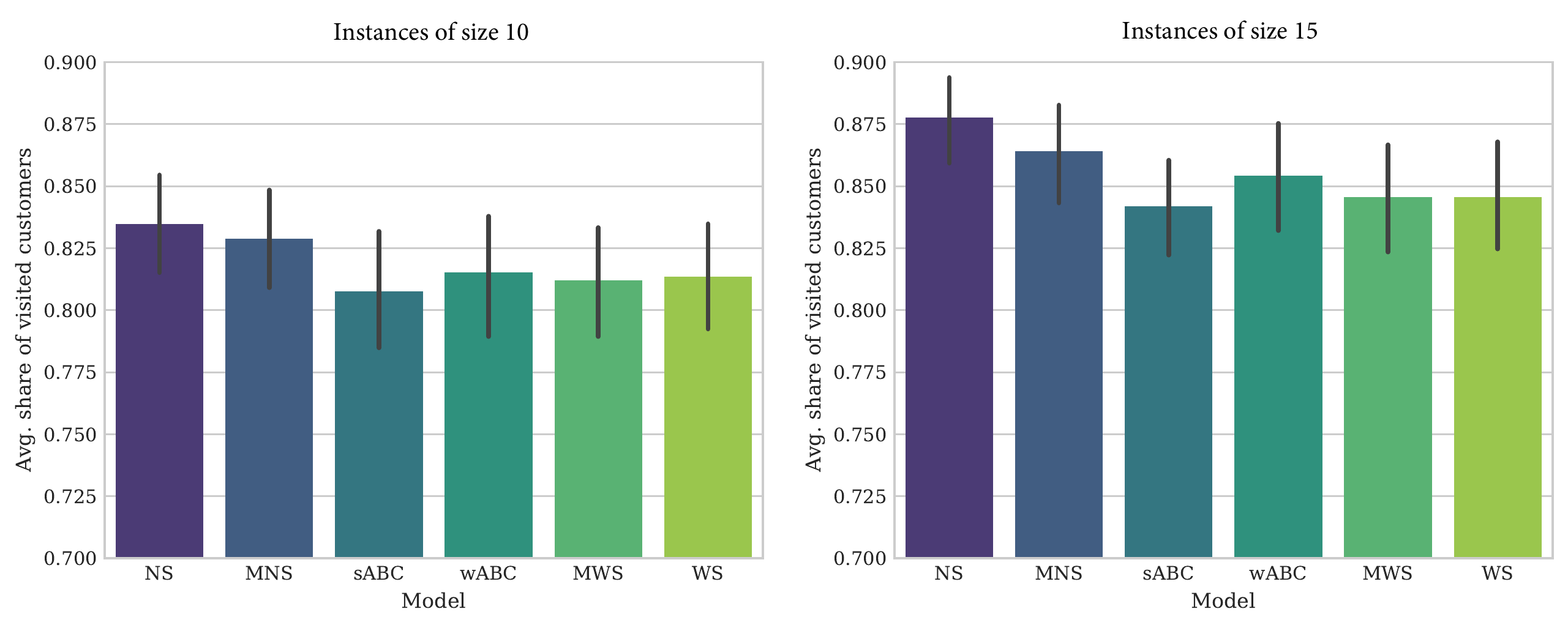}
  \caption{Average share of visited customers for instances with 10 (left) and 15 (right) customers.  The black line indicates 95~\% confidence interval. For a better readability the y-axis values start with values greater than zero.}
  \label{fig:SmallVisited}
\end{figure}

Figure~\ref{fig:SmallVisited} shows the share of visited customers averaged over the instances of size 10 and 15.
As expected, the number of visited customers is highest if the NS model is applied.
MNS delivers slightly worse values since the consideration of an unfavorably located mandatory customer in a solution prevents the service of two or more optional customers in some cases.
All other models aim to maximize the scores in different ways. Hence, in these solutions, fewer customers are visited in favor of serving far away customers with higher scores.
The sABC model is the most restrictive one in the sense that we enforce to first visit even the far away A and B customers before starting to ``fill-up'' the tours with B or C customers, who can be reached on the way.
Therefore, the sABC model often results in the fewest visited customers.


\paragraph{Realized scores:}

\begin{figure}[h!]
  \includegraphics[width=\textwidth]{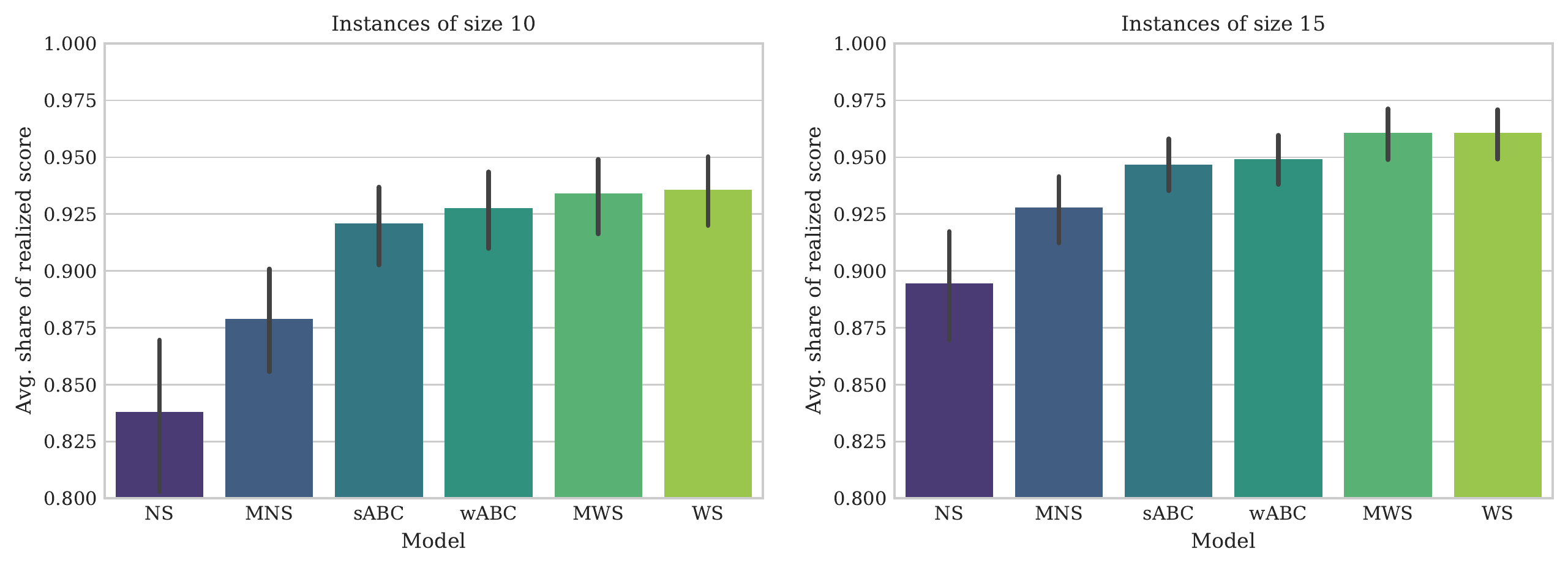} 
  \caption{Average share of realized score for instances with 10 (left) and 15 (right) customers. For a better readability the y-axis values start with values greater than zero.}
  \label{fig:SmallScores}
\end{figure}

Figure~\ref{fig:SmallScores} shows the share of realized scores for the same set of instances.
The NS model, which does not rely on any type of score information, results in the lowest share of realized score. 
For the other models, the graphs show a clear trend:
The more information expressing the profitability of the customer is used in a model, the better are the resulting shares of realized score.
Compared to the NS results, the shares can be considerably increased by adding promising customers to the set of mandatory customers in the MNS model.
On average, the models sABC and wABC yield further improvements and achieve results with similar shares of realized score.
As it can be expected, the WS model which maximizes the realized scores performs best.
The values of the MWS model are very close to the ones of WS.
There is just one instance in which a mandatory customers has lead to a lower realized score.

In summary, we can say that if only few information about customers is available, the manual selection of the most important customers as mandatory already improves the result considerably.
The models sABC and wABC also perform quite well considering the significantly smaller effort for the data acquisition in the prediction phase.

\paragraph{Varying score distributions:}

\begin{figure}[h!]
  \includegraphics[width=\textwidth]{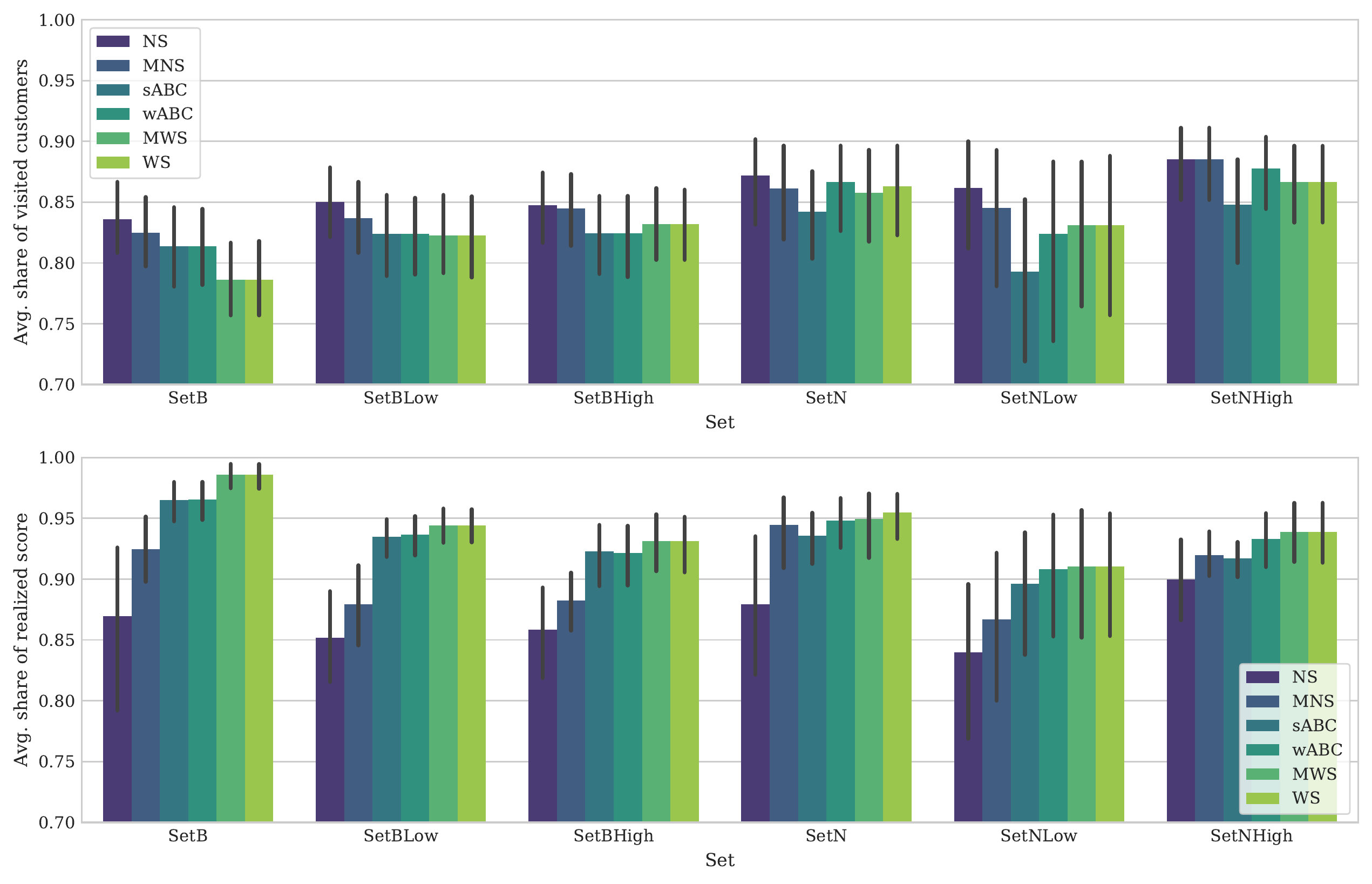} 
  \caption{Average share of visited customers (upper) and average share of realized score by set for instances with 10 and 15 customers. For a better readability the y-axis values start with values greater than zero.}
  \label{fig:SmallVisitedBySet}
\end{figure}

Figure~\ref{fig:SmallVisitedBySet} shows the share of visited customers on the upper graph and the share of realized scores differentiated by set and model.
The share of visited customers is different between the two base instance sets \Nesp and \Brill.
Within the sets, the different models behave consistently. 
Concerning the realization of scores, the {effect is consistent} over all score distributions:
To additionally consider information is positive; 
the more detailed the information is, the better are the results.
Furthermore, the results show that wABC consistently delivers good results over all scenarios.

\subsection{Results on large instance sets}

For analyzing the results of the 2MLS on the instances of the large instance sets, we only considered the best solution over the ten runs.

\paragraph{Impact of different models:}

\begin{figure}[h!]
  \includegraphics[width=\textwidth]{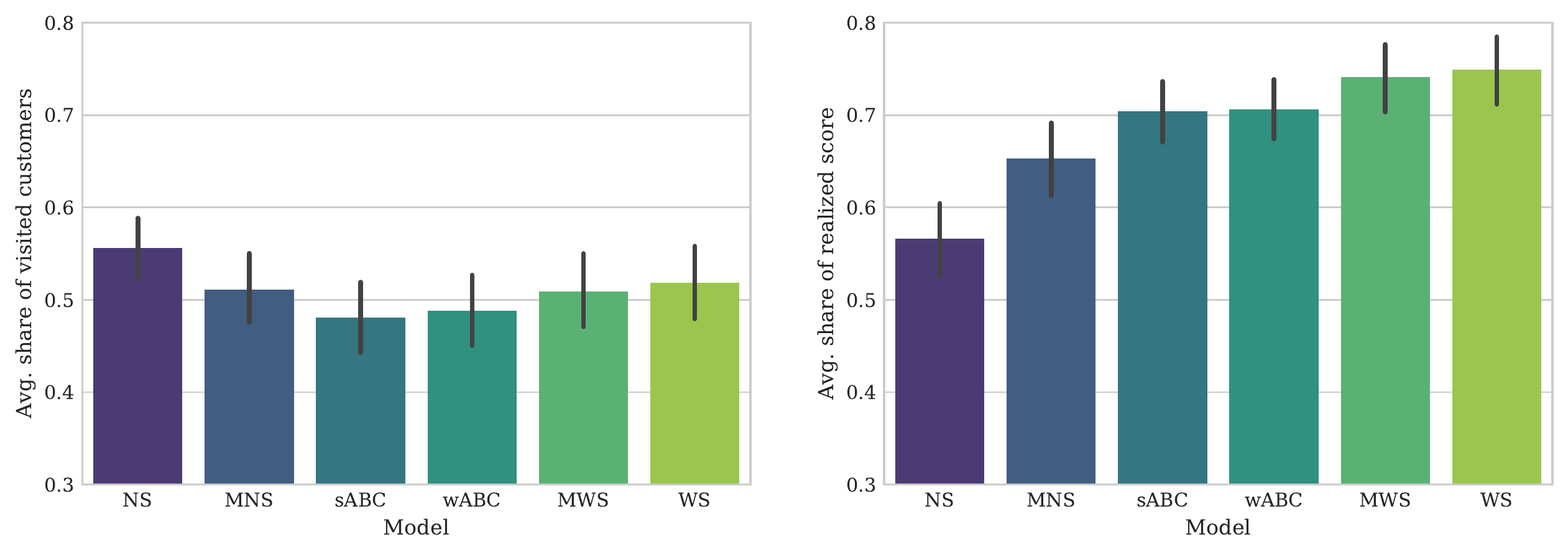} 
  \caption{Average share of visited customers (left) and average share of realized score (right) over all instances by model. For a better readability the y-axis values start with values greater than zero.}
  \label{fig:RelCustScoreAll}
\end{figure}

For a first overview, Figure~\ref{fig:RelCustScoreAll} shows the average share of visited customers, and the average share of realized score over all instances.
On this aggregated level it already becomes clear that the effects of different modeling variants identified on the small instance sets are confirmed by the heuristic solutions for the large instances.
Furthermore, on the small instances discussed above we obtained a comparatively high share of visited customers, with 80~\% to 95~\% (see Figure~\ref{fig:SmallVisited}).
On large instances, this share only lies around 50~\% (see Figure~\ref{fig:RelCustScoreAll}).
That means, the effects also seem to be independent of different levels of coverage.

\begin{figure}[h!]
  \includegraphics[width=\textwidth]{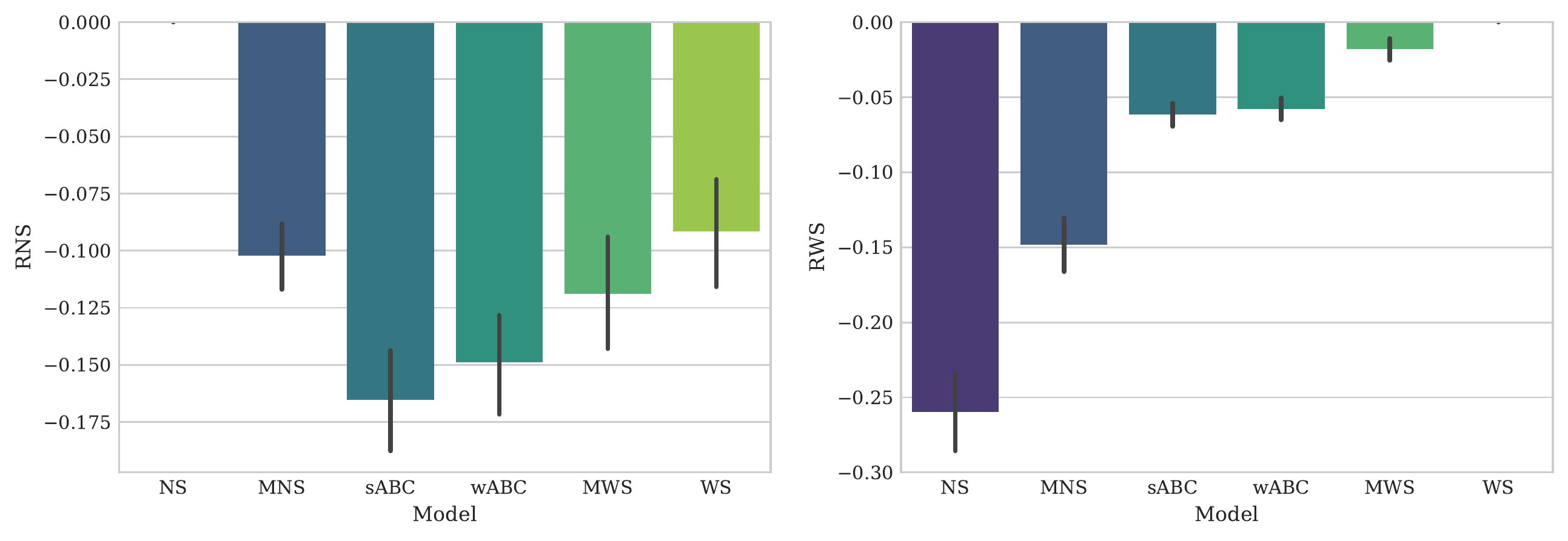} 
  \caption{Average RNS (left) and average RWS model (right)  over all instances.}
  \label{fig:RelNSRelWS}
\end{figure}
\begin{figure}[h!]
  \includegraphics[width=\textwidth]{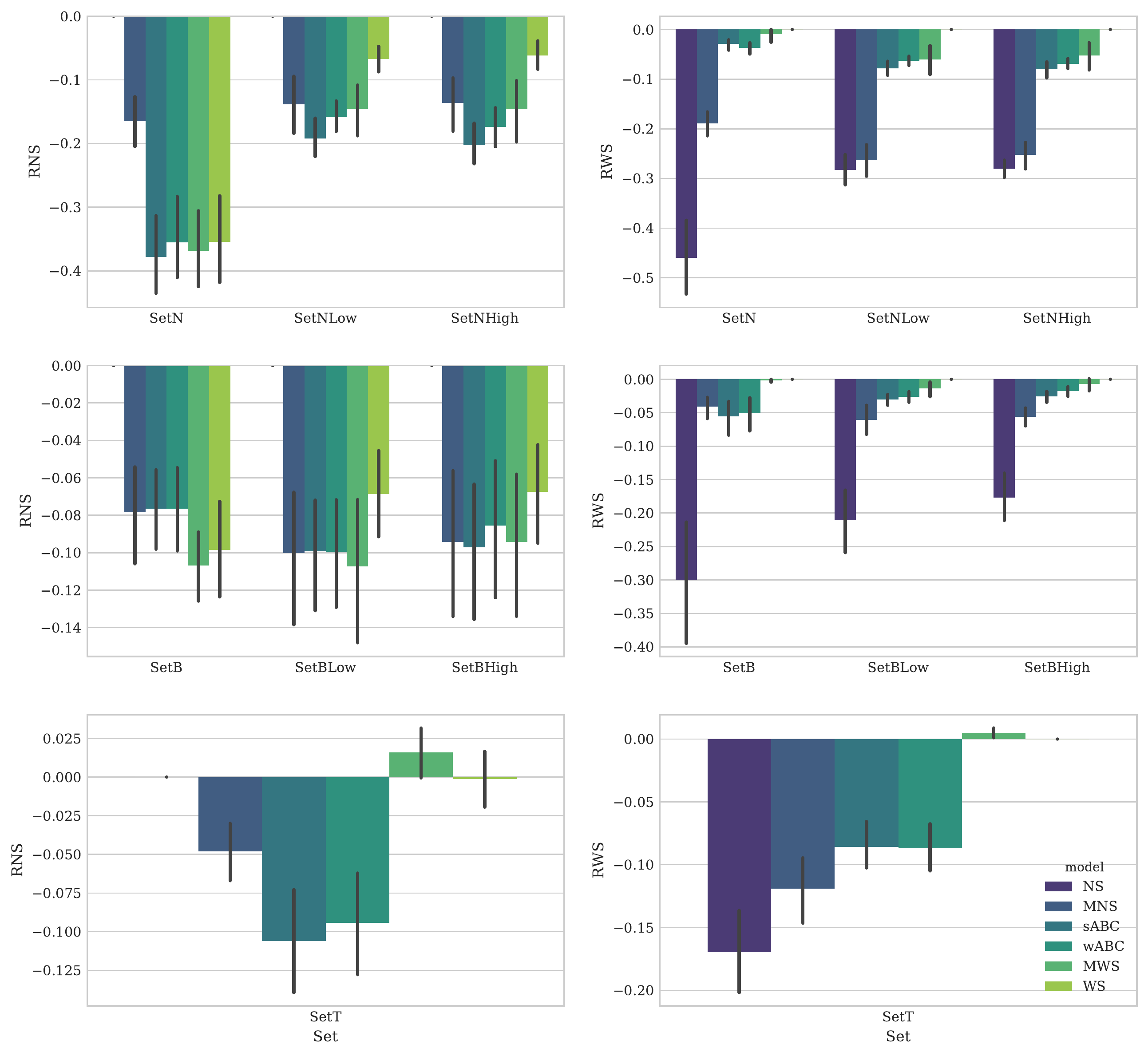} 
  \caption{Average RNS (left) and average RWS (right) by set and model.}
  \label{fig:RelNSRelWSBySet}
\end{figure}

A more detailed analysis is given in Figures~\ref{fig:RelNSRelWS} and~\ref{fig:RelNSRelWSBySet}.
To compare the relative performance of different models, we give the share of visited customers relative to the NS model, which directly maximizes this number.
This measure is indicated as ``RNS''.
Similarly, we relate the share of realized score to the WS model.
This measure is indicated as ``RWS''.
To take an example of Figure~\ref{fig:RelNSRelWS}, in MNS solutions 10~\% fewer customers are visited compared to NS solutions, while on average 15~\% less score is realized in MNS solutions compared to WS solutions.  


An interesting finding from these charts is that the sABC model is slightly better, i.e., it reaches higher realized scores than wABC, on instances with a strongly skewed distribution of scores, namely on sets \Nesp, \Brill, and \Tric.
On the sets with an equal distribution (\NespOne, \NespTwo, \BrillOne, \BrillTwo), the wABC model delivers slightly better results.
In our first, but false assumption, we attempted to explain this effect by the selection of the parameters of the wABC model:
For our computational study we used $p^A=15$, $p^B=5$, and $p^C=1$.
In \Nesp the average score over all instances is 929 for class A, 463 for class B, and 82 for class C. For \NespOne we have an average score of 833 for class A, 498 for class B, and 159 for class C. 
This means, in \Nesp an estimator for $p^A$ would be $929 / 82 = 11.3$, for $p^B=463 / 82 = 5.6$, and $p^C=1$. 
For \NespOne the parameters could be set to $p^A=5.2$, $p^B=3.1$, and $p^C=1$.
However, the selected parameters for our computational study were considerably closer to the estimators of \Nesp, for which the sABC model performs better than the wABC variant.

A more likely explanation for the effect is that in \Nesp, \Brill, and \Tric the number of A and B customers is considerably smaller than in the variants with equal score distributions.
{This means that the probability of leaving out a C customer instead of a B customer in favor of A customers is much higher than in the other investigated scenarios.}
Furthermore, unvisited C customers have a lower average score. On average, C customers have a score of 82 in \Nesp, 152 in \NespOne, and 4264 in \NespTwo.
In other words, in case of a strongly skewed score distribution with a low share of very important customers in classes A and B, the strictly hierarchical objective function of sABC seems a good choice.
If the scores are rather equally distributed, the selection of wABC seems better, even if from this computational study no recommendations concerning a good parametrization for wABC can be made. 

In general, we can state that the sABC and the wABC model both deliver good results considering the substantially smaller effort on the data acquisition side.
Furthermore, if no data is available, it is worth to mark the important customers as mandatory manually:
{The resulting scores with the MNS model are considerably better compared to the NS model.}
As can be expected, this effect is especially high in the case of strongly skewed score distributions.

\paragraph{Interaction of model and 2MLS:}

Two results indicate, in some cases, the combination of the applied model and the 2MLS can have a positive impact on the results:
In some cases in \Tric, the number of visited customers is increased by using MWS or WS compared to the NS case, which represents the lower bound in case of an exact solution method.
Also in some cases of \Tric, the score of the MWS model is higher than the one of the lower bound model WS.
There are two possible explanations:
First, mandatory customers in a MWS model, who correspond to the most important customers for reaching high score solutions, seem to produce a favorable structure.
This applies to the multi-start procedure (as seed candidates) and the insertion phase, in which they are inserted before optional customers.
Second, considering scores in the MWS and WS models leads to more diversified solutions, and higher diversification of solutions usually prevents the search from remaining in locally optimal solutions.
However, none of the effects can be consistently observed over the different instance sets.
Hence, none of the models is preferable concerning the resulting quality of the 2MLS.

\subsection{Sensitivity to prediction and classification errors}

\begin{figure}[tb!]
  \centering
  \includegraphics[width=.6\textwidth]{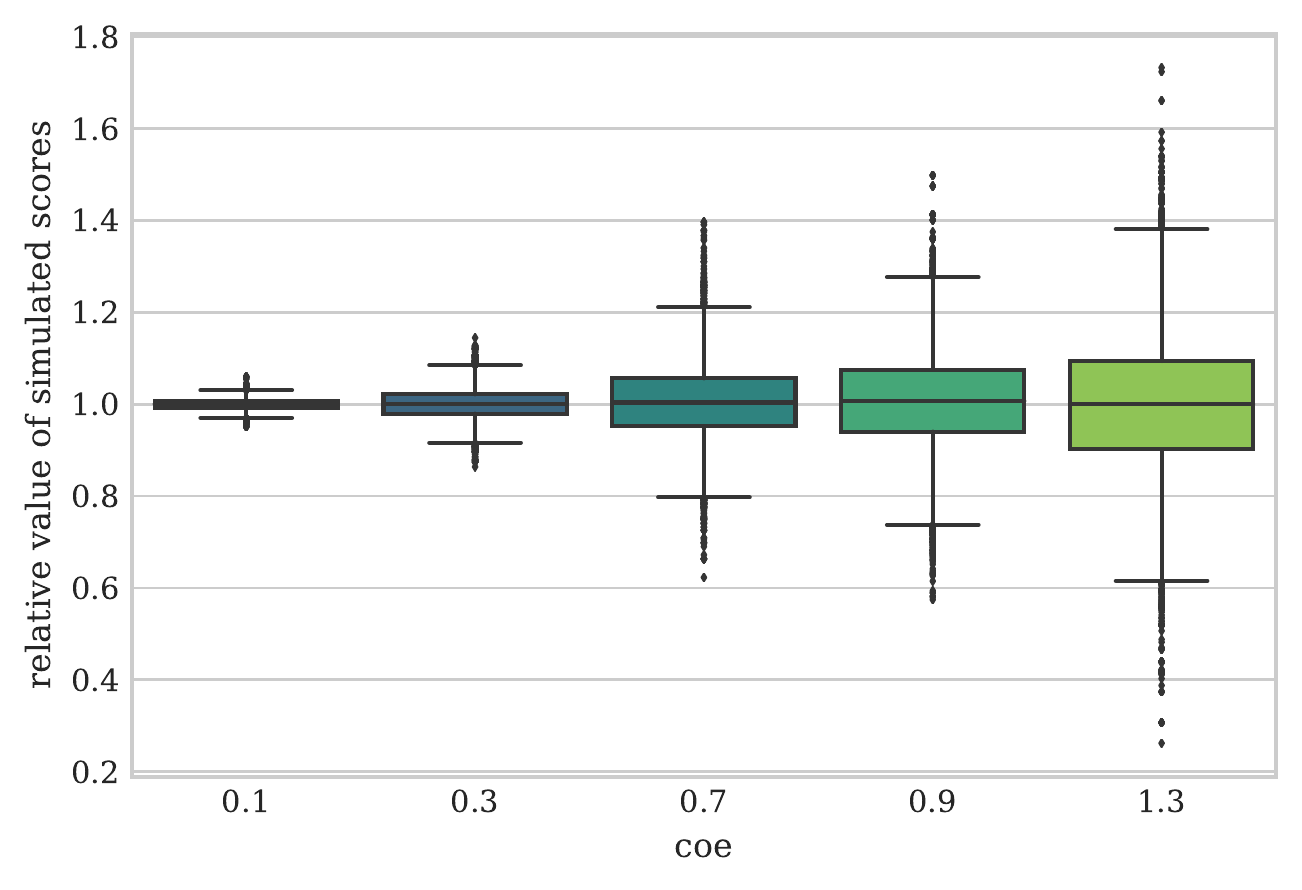} 
  \caption{Simulated score realizations relative to predicted outcome for WS models.}  
  \label{fig:resultsSensitivityWS}                                        
\end{figure}

In our sensitivity analysis, we study the performance of the different modeling approaches given different simulated score realizations $p_i^{sim}$. 
Figure \ref{fig:resultsSensitivityWS} represents the value of the simulated realizations relative to the scores used in the planning phase for the WS model, i.e., it measures the prediction error associated with a plan assuming that the actual value of customers is revealed after a visit (see Section \ref{subsec:pred_quality}).
As expected from the experiment setup, the resulting values are symmetrically distributed around 1 and vary significantly for larger values of the $coe$.
 
{To compare different models considering the robustness towards errors, we define a measure $\text{RWS}^{sim}$, which gives the sum of the simulated score realization for any one model as a percentage of the result achieved with the WS model.
Figure \ref{fig:resultsSensitivityBase} gives the results for \Nesp and \Brill, i.e., for the base instances with skewed score distributions.
For both sets, in case of high quality predictions, i.e., in case of small values of $coe$, the advantage of explicitly considering the predicted score information in the WS objective value is obvious and in line with our prior observations. 
For nearly all test instances, all other models consistently result in worse scores. }

{In case of larger prediction errors, however, the $\text{RWS}^{sim}$ measure for sABC and wABC models on average come pretty close to 1, i.e., the achieved outcome is of a similar quality as the one obtained using WS.
Furthermore, it becomes more likely that the WS solutions are outperformed.
In case of \Brill, even the very simple MNS model comes pretty close to 1 and performs surprisingly good in case of $coe=1.3$. 
{This shows that it is no longer beneficial to take more information into account when prediction quality is low. 
Especially in cases of highly skewed score distributions, relatively simple classification-based methods provide very robust results as the most important customers can still be reliably identified.
}


\begin{figure}[tb!]
  \includegraphics[width=\textwidth]{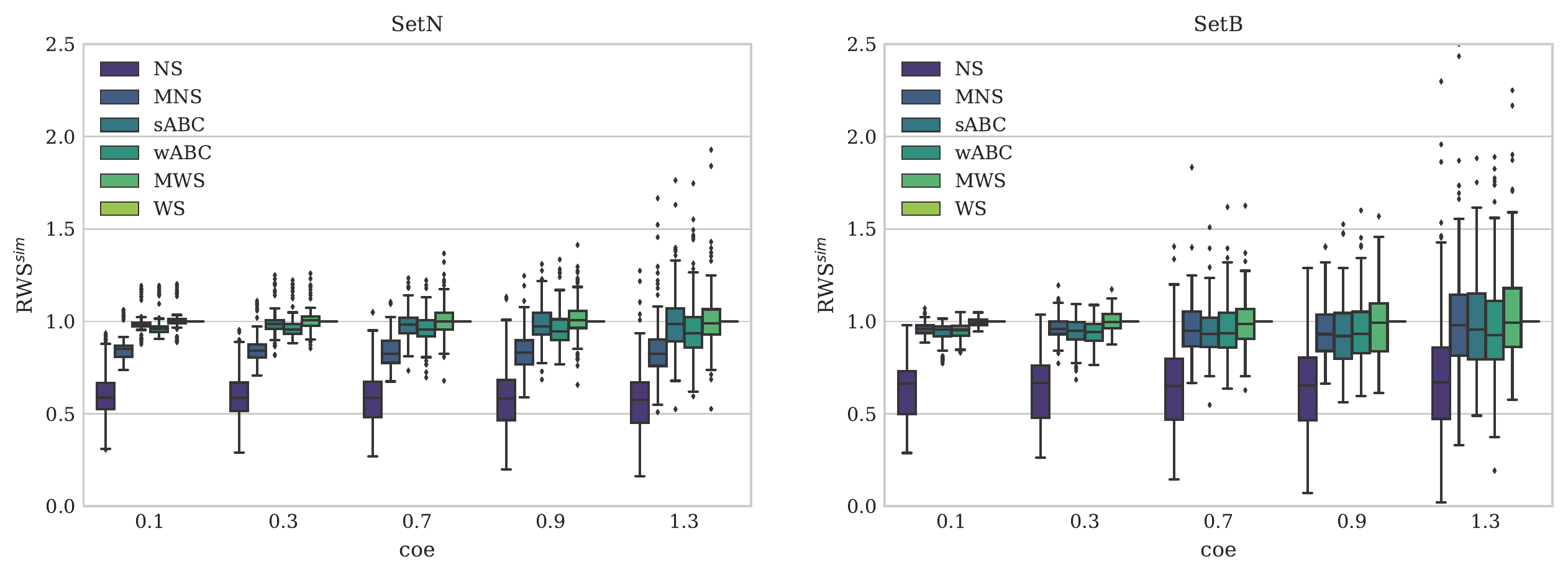} 
  \caption{Impact of prediction errors depending on the applied model for \Nesp and \Brill.}  
  \label{fig:resultsSensitivityBase}                                        
\end{figure}

{
The results on \NespOne and \BrillOne are given in Figure \ref{fig:resultsSensitivityLow}.
Note that \NespTwo and \BrillTwo, as similarly homogeneous score distributions, exhibit the same effects.
Again, we can observe that WS reliably achieves the best results for low $coe$.
Compared to \Nesp, the average $\text{RWS}^{sim}$ of sABC, wABC and MWS models is less than 1, indicating that WS is beneficial even in case of larger prediction errors.
In contrast, in \BrillOne, $\text{RWS}^{sim}$ quickly increases with increasing $coe$ for a larger number of scenarios. 
}

\begin{figure}[tb!]
  \includegraphics[width=\textwidth]{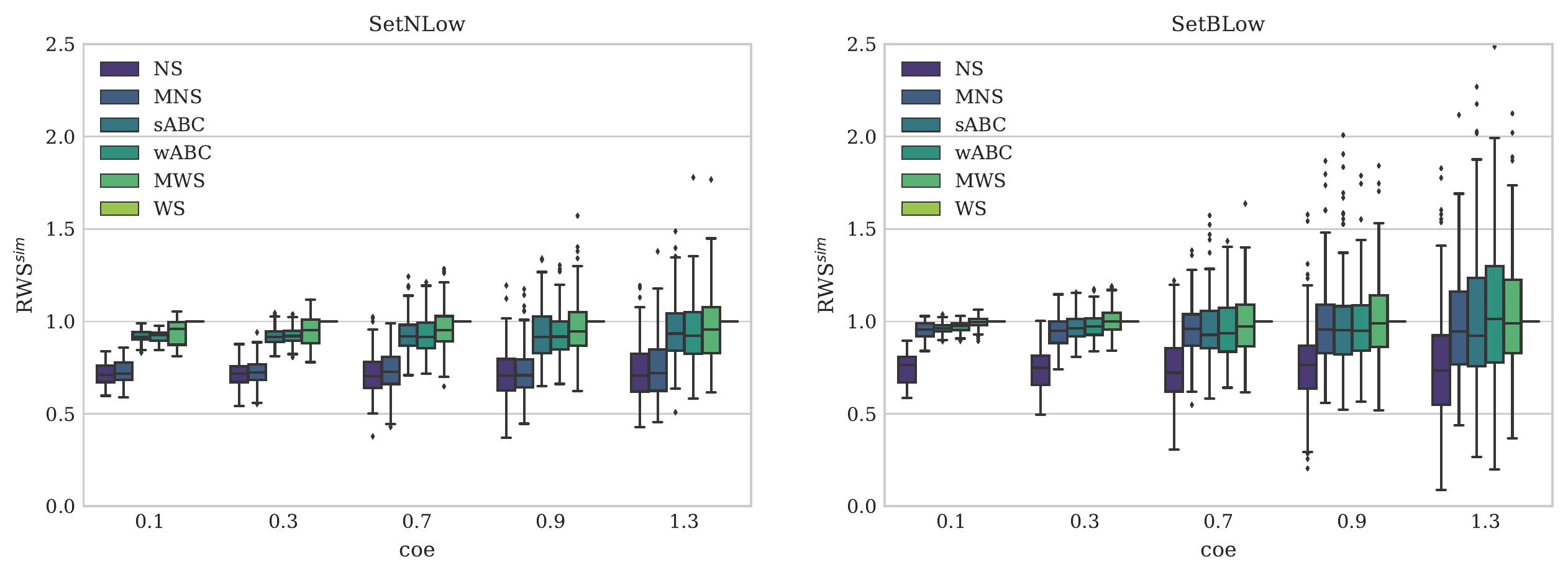} 
  \caption{Impact of prediction errors depending on the applied model for \NespOne and \BrillOne.}  
  \label{fig:resultsSensitivityLow}                                        
\end{figure}

\section{Managerial insights and conclusion}\label{sec:outlook}

As shown in the preceding discussions, there exist scoring methods for customer selection ranging from judgemental scoring and classification of customers conducted ``manually'' to the prediction of customer response functions for concrete offers based on an extensive database using machine learning.
The resulting scores vary in the level of aggregation and the content of information they carry.
However, in contrast to direct marketing campaigns by mail or phone calls, in the case of sales representatives working in the field, the final customer selection is strongly intertwined with tour planning decisions. 
Since the scoring methods also differ considerably in the required data basis and the effort for the data acquisition and maintenance, it is worth to consider the results of both stages to determine the necessary level of detail of the customer scoring.
{This is especially true if there is no advanced customer relationship management system in place that tracks the sales process on a very detailed level, and if no team dedicated to maintaining data quality has been established.}

{To study this interplay between customer scoring and tour planning, we proposed several modeling variants of the multi-period orienteering problem that each incorporate a specific category of customer scoring method.}
Based on the results of an extensive computational study {using several real-world test instances}, the following insights might be helpful for practitioners:

\begin{itemize}
	\item If the predicted scores represent the response of customers to an offer (e.g. additional sales) and if they are reasonably accurate, the effort in data preparation for using the WS model is justified by a consistent outperformance over all other models. In our use cases, the advantage was observable until a $coe$ of $0.7$.
	If the prediction errors of the response functions increase, classification and assessment based scoring methods and tour planning variants offer robust outcomes with less need for {high resolution,} accurate data.
	\item If the available data basis is smaller or contains no detailed history and only allows for scoring methods of category rank or class, both sABC, and wABC models deliver very good results.
	This is especially the case, if the scores are strongly skewed, i.e. if there are few customers with comparatively high scores, and many others with low scores, as it was the case in our real-world instances.
	If the score distribution is strongly skewed, it is more important to identify the important customers than differentiating the less important candidates, e.g. prospects.
	\item If the data basis is very small, e.g. in a new district, it is worth to incur great expense to identify the most important customers -- even manually -- and select the MNS model variant. 
	Compared to considering no scores at all in the NS model, MNS improves the share of realized scores in many situations considerably while visiting significantly fewer customers.
\end{itemize}

In general, we can state that in many situations customer classes seem to provide enough guidance for reaching good results in the tour planning stage.
Even the identification of a few very important customers guides the search surprisingly well.
For tactical and strategic decisions, usually, predictions on an even more aggregated level are used. 
This shows that for managing an efficient data acquisition and maintenance process it is highly recommended to consider both stages of the short term planning, the customer scoring (or prediction) and the tour planning (or decision) stage.


It is also worth mentioning that the consideration of mandatory customers can be used for different purposes:
Private appointments of the sales representative or arranged dates with customers, as well as customers who must be visited regularly can be modeled as mandatory customers and can be fixed to a day.
From our experience, the possibility of considering private dates of sales representatives is an absolutely crucial feature for the acceptance of tour planning support tools within the sales force.

In future work, we will consider tour planning approaches explicitly considering the stochastic behavior by applying modeling variants from stochastic programming or robust optimization.
This extends the point estimations for the customer profit applied in this work.

\begin{appendix}

\pagebreak

\section{Appendix}

\subsection{Main components of  2MLS}
\label{sec:2mls}

The 2MLS of \cite{Glock2020} relies on the concepts of the adaptive large neighborhood search (ALNS) first introduced in \cite{Ropke2006} and combines it with a multi-start phase for quickly identifying promising starting solutions. 
For solving the profitable sales representatives tour problem, we left out all strategies related to the multi-vehicle case and used the same parameter settings, if applicable.

Similar to the solution concepts discussed in Section \ref{sec:orienteering_literature}, we combine neighborhoods that focus on improving the sum of scores of the selected customers with others that seek to improve the efficiency of the obtained tours.
These are implemented in form of adapted removal and insertion strategies in the ALNS, which are briefly discussed below. 
We then introduce the multi-start strategy and the reheating scheme that guides the search.

\paragraph{Removal strategies}

Removal strategies remove a given percentage of customer nodes from a solution following relatively simple heuristic rules.
In our implementation, unvisited customers have to be explicitly removed from a solution as well.
Otherwise, they cannot be re-inserted in the subsequent repair step.
Unless otherwise noted, the unvisited customers to be removed are neighbors of those removed from the tours.
The purpose of this is to limit the size of the neighborhood in each iteration and to reduce the computational effort due to evaluating customers that are not likely to be re-inserted. 
Our 2MLS uses six different removal strategies for solving the MPOP:
\begin{description}
\item[RND-NN] removes visited customers at random.
\item[SEQU-NN] removes sequences of customers from the tours.
\item[SCORE-DELTA] removes the planned customers with the lowest scores as well as the unvisited customers with the highest scores from a solution.
\item[SKEL-TOUR] randomly selects one tour and removes all customer visits within this tour in segments of a fixed length, leaving one customer in between subsequent segments. This results in a ``skeleton'' tour with the same basic structure as the initial one.
\item[WORST-DETOUR] removes customers associated with high travel times.
\item[WORST-ANGLE] selects customers that come with sharp angles in a tour as well as the preceding and succeeding customers.
\end{description}

\paragraph{Insertion strategies}

Insertion strategies iteratively select one of the removed customers and try to reinsert it into a tour.
Strategies vary in the order of insertions. 
If there are mandatory customers, i.e., $\mathcal{N}^{M} \neq \emptyset$, they are inserted first; if the insertion is not feasible, the repair phase ends.
The insertion strategies are:
\begin{description}
\item[RND] inserts customers in random order.
\item[MAX-SCORE] inserts the remaining customer with the highest score first.
\item[SCORE-RATIO] calculates the ratios of the customers' scores with the associated detours for the best insertion positions and inserts the one with the best ratio first.
\item[SCORE-RATIO2] similarly computes the score ration but uses the squared score values.
\item[GREEDY] inserts customers in increasing order of the necessary detour, irrespective of their score.
\end{description}
All strategies (except RND) are slightly randomized so that in each iteration, an arbitrary customer is selected for insertion with a probability of $0.05$ to encourage diversification.

The selected customers are inserted at the position that is currently associated with the lowest detour with respect to the emerging tours.
If no position remains that is feasible with respect to the working time constraint, the customer is left unvisited in the solution.

\paragraph{Seeds and multi-start}

Initial solutions are created in a multi-start approach.
First, we select seed nodes for each day from the set of customers using a k-means++ variant, which means that the likelihood of a customer to be selected as seed increases with its distance towards other seeds and the home location \citep{Arthur2006}. 
If there are mandatory customers, the seeds are chosen from the set $\mathcal{N}^{M}$.
Each seed initializes one tour.
Then, one of the insertion strategies is applied to build tours for each day.
This is repeated 10 times using different insertion strategies.
The best solution obtained this way is used as a starting solution for the ALNS.

\paragraph{Acceptance and stopping criterion}
During the search, improving solutions are always accepted.
Similar to \cite{Vidal2015}, we use travel distance as a tiebreaker for solutions with same scores.
This way, the search prefers solutions visiting the same customers at a lower cost.
We furthermore apply a reheating scheme related to simulated annealing to guide the search.
This scheme accepts solutions that come with a decrease in solution quality with a probability that increases with the number of iterations since the best solution was last improved.
Decreases in quality are only accepted if the search has likely resulted in a local minimum that cannot be left otherwise.
The search stops when no improvements over the current best solution have been found in 300 iterations.

%

\end{appendix}

\clearpage
\bibliographystyle{agsm}
\bibliography{lit}

\end{document}